\documentclass[a4paper,12pt,reqno]{amsart}
\usepackage{amssymb}
\usepackage[dvips]{graphicx}
\usepackage{hyperref}

\input epsf.tex
\newdimen\xsize
\newdimen\oldbaselineskip
\newdimen\oldlineskiplimit
\def\restorelineskip{\baselineskip=\oldbaselineskip%
\lineskiplimit=\oldlineskiplimit}
\def\putm[#1][#2]#3{
\hbox{\vbox to 0pt{\parindent=0pt%
\vskip#2\xsize\hbox to0pt{\hskip#1\xsize $#3$\hss}\vss}}}%
\long\def\Line#1{\hbox to \hsize{#1}}
\def\putt[#1][#2]#3{
\vbox to 0pt{\noindent\hskip#1\xsize\lower#2\xsize%
\vtop{\restorelineskip#3}\vss}}

\makeatletter
\def\xbig[#1]#2{{\hbox{$\m@th\left#2\vbox to#1\xsize{}%
\right.\n@space$}}}
\makeatother
\def\xlar[#1]#2{%
\smash{\mathop{ \hbox to #1\xsize{\leftarrowfill}}\limits^{#2}}}
\def\xrar[#1]#2{%
\smash{\mathop{ \hbox to #1\xsize{\rightarrowfill}}\limits^{#2}}}
\def\xline[#1]{\hbox to #1\xsize{\leaders\hrule\hfill}}

\DeclareFontFamily{U}{rsf}{\skewchar\font'177}%
\DeclareFontShape{U}{rsf}{m}{n}{<-6>rsfs5<6-8>rsfs7<8->rsfs10}{}%
\DeclareFontShape{U}{rsf}{b}{n}{<-6>rsfs5<6-8>rsfs7<8->rsfs10}{}%
\DeclareMathAlphabet\RSFS{U}{rsf}{m}{n}
\SetMathAlphabet\RSFS{bold}{U}{rsf}{b}{n}
\def\sf#1{{\mathsf{#1}}}

\def\slsf{\slshape \sffamily }

\def\msmall#1{\mathchoice{\hbox{\small$\displaystyle {#1}$}}{#1}{#1}{#1}}

\def\bb{{\mathbb B}}

\def\cc{{\mathbb C}}

\def\rr{{\mathbb R}}

\def\nn{{\mathbb N}}

\def\pp{{\mathbb P}}
  
\def\zz{{\mathbb Z}}

\def\adyn{\sf{1}}

\def\c{\sf{c}}
\def\C{\sf{C}}

\def\dim{\sf{dim}\,}

\def\ind{\sf{ind}}
\def\j{\sf{j}}

\def\exp{\sf{exp}}

\def\grad{\sf{grad}}

\def\id{\sf{Id}}

\def\inter{\sf{Int}}

\def\lim{\mathop{\sf{lim}}}

\def\log{\sf{log}\,}

\def\max{\sf{max}}

\def\Reg{\sf{Reg}\,}

\def\Sing{\sf{Sing}\,}

\def\tr{\sf{tr}\,}

\def\v{{\mathrm{v}}}
\def\w{{\mathrm{w}}}

\def\z{\sf{z}}

\def\eps{\varepsilon}

\def\<{\langle}\let\la=\<
\def\>{\rangle}\let\ra=\>

\def\d{\partial}
\def\dbar{{\barr\partial}}

\def\ddef{\mathrel{{=}\raise0.3pt\hbox{:}}}
\def\deff{\mathrel{\raise0.3pt\hbox{\rm:}{=}}}

\def\fraction#1/#2{\mathchoice{{\msmall{ #1\over#2}}}%
{{ #1\over #2 }}{{#1/#2}}{{#1/#2}}}
\def\norm#1{\left\Vert{#1}\right\Vert}

\def\le{\leqslant}

\def\longpoints{\leaders\hbox to 0.5em{\hss.\hss}\hfill \hskip0pt}
\def\stateskip{\smallskip}
\def\state#1. {\stateskip\noindent{\bf#1. }} 
\def\statep#1. {\stateskip\noindent{\bf#1 }} 
\def\proof{\state Proof. \2}

\def\Chi{\raise 2pt\hbox{$\chi$}}

\def\ie{\hskip1pt plus1pt{\sl i.e.\/,\ \hskip1pt plus1pt}}

\def\sli{{\sl i)} } 
\def\slii{{\sl i$\!$i)} } 
\def\sliii{{\sl i$\!$i$\!$i)} }

\def\barr#1{\mskip1mu\overline{\mskip-1mu{#1}\mskip-1mu}\mskip1mu}
\def\Chi{\raise 2pt\hbox{$\chi$}}

\let\phI=\phi\let\phi=\varphi\let\varphi=\phI
\let\cal=\mathcal
\def\cala{{\cal A}}

\def\cale{{\cal E}}
\def\calf{{\cal F}}

\def\calh{{\cal H}}

\def\calk{{\cal K}}

\def\calm{{\cal M}}
\def\caln{{\cal N}}
\def\calo{{\cal O}}

\def\eps{\varepsilon}

\def\d{\partial}
\def\dbar{{\barr\partial}}
\def\1{{1\mkern-5mu{\rom l}}}

\def\ge{\geqslant}

\def\fraction#1/#2{\mathchoice{{\msmall{ #1\over#2}}}%
{{ #1\over #2 }}{{#1/#2}}{{#1/#2}}}

\def\le{\leqslant}

\newcommand{\2}{\thinspace}

\def\qed{\ \ \hfill\hbox to .1pt{}\hfill\hbox to .1pt{}\hfill $\square$\par}

\def\comment#1\endcomment{}

 \belowdisplayskip=\abovedisplayskip

\def\lineeqqno(#1){\hfill\llap{\vbox to 10pt%
{\vss\begin{align} \eqqno(#1)\end{align}\vss}}\vskip1pt}

\advance\headheight 1.2pt

\def\ShowwLLabel#1{}

\def\thechpt{\Roman{chpt}}
\def\newchapt[#1]#2{%
\refstepcounter{chpt}\setcounter{subsection}{0}%
\setcounter{thm}{0}\setcounter{defi}{0}%
\setcounter{rema}{0}\setcounter{exrc}{0}%
\renewcommand{\thesubsection}{\thechpt.\arabic{subsection}}%
\newpage
\section*{\begin{center}\huge \bf Chapter \thechpt\\
#2 \end{center}}\label{#1}%
\ \smallskip%
\addcontentsline{toc}{chpt}{Chapter \thechpt. #2}%
\markboth{Chapter \thepart}{#2}%
}
\def\newsect[#1]#2{\refstepcounter{section}\setcounter{equation}{0}%
\renewcommand{\thesubsection}{\arabic{section}.\arabic{subsection}}%
\section*{\arabic{section}.
#2}\vspace{-20pt}\label{#1}\vspace{20pt}%
\markboth{Section \arabic{section}}{#2}}

\def\newlect[#1]#2{\refstepcounter{section}%
\renewcommand{\thesubsection}{\arabic{section}.\arabic{subsection}}%
\section*{Lecture \arabic{section}\\
#2}\label{#1}%
\markboth{Lecture \arabic{section}}{#2}}

\def\newprg[#1]#2{\refstepcounter{subsection}%
\subsection*{{\thesubsection.\ #2}} \label{#1}%
}

\def\newappx[#1]#2{%
\refstepcounter{appx}\setcounter{section}{0}%
\renewcommand{\thesubsection}{A\arabic{appx}.\arabic{subsection}}%
\section*{Appendix \arabic{appx}\\ #2}
\label{#1}%
\markboth{Appendix A\arabic{appx}}{#2}
}

\newtheorem{thm}{Theorem}[section]
   \def\newthm#1{\begin{thm}\label{#1}}

\newtheorem{nnthm}{Theorem}
   \def\newthm#1{\begin{nnthm}\label{#1}}

\newtheorem{lem}{Lemma}[section]
   \def\newlemma#1{\begin{lem} \label{#1}}

\newtheorem{prop}{Proposition}
   \def\newprop#1{\begin{prop}\label{#1}}

\newtheorem{nnprop}{Proposition}[section]
   \def\newprop#1{\begin{nnprop}\label{#1}}

\newtheorem{corol}{Corollary}
   \def\newcorol#1{\begin{corol} \label{#1}}

\newtheorem{nncorol}{Corollary}[section]
   \def\newcorol#1{\begin{nncorol} \label{#1}}

\newtheorem{defi}{Definition}[section]
   \def\newdefi#1{\begin{defi} \label{#1}\rm }

\newtheorem{exmp}{Example}[section]
   \def\newexmp#1{\begin{exmp} \label{#1}\rm }

\newtheorem{nnexmp}{Example}
   \def\newexmp#1{\begin{nnexmp} \label{#1}\rm }

\newtheorem{exrc}{Exercise}
   \def\newexrc#1{\begin{exrc} \label{#1}\rm }

\newtheorem{rema}{Remark}[section]
   \def\newrema#1{\begin{rema} \label{#1}\rm }

\newtheorem{nnrema}{Remark}
   \def\newthm#1{\begin{nnrema}\label{#1}}

\def\eqqno(#1){\label{(#1)}}
\def\eqqref(#1){(\ref{(#1)})}

\title{Bochner-Hartogs type extension theorem for roots\\
and logarithms of holomorphic line bundles}
\author{S. Ivashkovich}
\date{\today}

\address{
Universit\'e de Lille-1, UFR de Math\'ematiques, 59655 Villeneuve
d'Ascq, France} \email{ivachkov@math.univ-lille1.fr}
\address{IAPMM Nat. Acad. Sci. Ukraine
Lviv, Naukova 3b, 79601 Ukraine}


\begin{document}
\begin{abstract}
We prove an extension theorem for roots and logarithms of
holomorphic line bundles across strictly pseudoconcave boundaries:
they extend in all cases except one, when dimension and Morse index
of a critical point is two. In that case we give an explicit
description of obstructions to the extension.
\end{abstract}

\maketitle

\setcounter{tocdepth}{1}
\tableofcontents

\newsect[sect.INT]{Introduction.}

\newprg[prgINT.stat]{Bochner-Hartogs type extension}
Let $X$ be a  complex manifold of dimension $n\ge 2$, $U$ some
domain in $X$ and let $\rho: U \to \rr$ be a smooth, strictly
plurisubharmonic, Morse function on $U$, such that for some $t_0$
the level set $\Sigma_{t_0} \deff \{\rho = t_0\}\subset U$ is
connected and compact. In the sequel, for $t$ close to $t_0$, we
denote  by $U_t^{+} \deff \{x\in U: \rho (x) >t\}$ the upper level
set of $\rho$ and by  $U_t^{-}\deff \{x\in U: \rho (x)<t\}$ - the
lower level set. Let a holomorphic line bundle $\caln$  on $U$ be
given. We say that $\caln$ admits a $k$-th root on $U_{t_0}^{+}$ if
there exists a holomorphic line bundle $\calf$ on $U_{t_0}^{+}$ such
that $\calf^{\otimes k} \cong \caln|_{U_{t_0}^{+}}$. The main goal
of this paper is to understand whether or not a root $\calf$ can be
extended from $U_{t_0}^{+}$ to a neighborhood of $\Sigma_{t_0}$. It
occurs that the answer depends on dimension $n$ of $X$ and on index
$\ind\rho|_{\c}$ of a critical point $\c$ of $\rho$ (if such exists
on $\Sigma_{t_0}$). It is clear that without loss of generality we
can suppose that the level set $\Sigma_{t_0}$ is either smooth (\ie
$\grad \rho|_{\Sigma_{t_0}}\not=0$) or, $\rho$ has exactly one
critical point $\c$ on $\Sigma_{t_0}$ of index (necessarily) $0\le
\ind\rho|_{\c}\le n=\dim_{\cc} X$. Also $U$ should be viewed as a
neighborhood of $\Sigma_{t_0}$. Our first result can be stated as
follows.

\begin{nnthm}
\label{root-ext-thm}
Suppose that $\Sigma_{t_0}$ is  either smooth

\smallskip\sli or, $n\ge 3$ and $\ind\rho|_{\c}$ is arbitrary;

\smallskip\slii or, $n=2$ and $\ind\rho|_{\c}\not=2$.

\smallskip Then $\calf$ extends to a holomorphic line bundle to a
neighborhood of $\Sigma_{t_0}$ and stays there to be a $k$-th root
of $\caln$.
\end{nnthm}

Now let us turn to the exceptional case $n=2$ and $\ind = 2$.
Consider the following model situation. Denote by $z_j = x_j+iy_j,
j=1,2$ the coordinates in $\cc^2\cong \rr^4$ and let $K^2 = \{|x_j|,
|y_j|<1: j=1,2\}$ be the unit cube in $\cc^2$. Set  $D\deff
K^2\setminus \{x_1=x_2=0\}$.  For every $k\ge 2$ and $1\le l <k$ we
are going to construct a holomorphic line bundle $\calf_k(l)$ on $D$
such that $\calf_k(l)^{\otimes k}\cong \calo$, \ie $\calf_k(l)$ is a
$k$-th root of the trivial bundle, which doesn't extend to $K^2$.

\smallskip Let $C\deff \{(x,y): y_1=y_2=0,  x_1^2+x_2^2=1/2\}$ be
the generator of $H_1(D,\zz)$. Since the totally real plane
$\{x_1=x_2=0\}$ is removable for holomorphic functions in $D$ and
$K^2$ is simply connected the exponential $\exp : H^0(D,\calo)\to
H^0(D,\calo^*)$ is surjective. Moreover, $H^2(D,\zz)=0$ and
therefore the sequence
\begin{equation}
\eqqno(exact0)
0 \to H^1(D,\zz)  \buildrel \j \over \to H^1(D,\calo) \buildrel \exp
\over \to H^1(D,\calo^*)  \buildrel c_1
\over \to 0
\end{equation}
is exact. Here $\j : H^1(D,\zz)  \to H^1(D,\calo)$ is the natural
imbedding. Denote by $A$ the cohomology class in $H^1(D,\zz)$ such
that $<A,C>=1$, and by $B$ denote its image $\j (A)$ in
$H^1(D,\calo)$. Set further $F_k(l) \deff \exp(\frac{l}{k} B)\deff
e^{\frac{2\pi i l}{k} B}\in H^1(D,\calo^*)$ and denote by
$\calf_k(l)$ the corresponding holomorphic line bundle. Since
$(F_k(l))^k = (\exp\circ\j) (A) =0$ we have that
$\calf_k(l)^{\otimes k} \cong \calo$, \ie $\calf_k(l)$ is a $k$-th
root of the trivial bundle on $D$. At the same time $F_k(l)$ is not
zero itself, because $l/k\cdot B$ cannot be an image of anything
from $H^1(D,\zz)$, and as such the bundle $\calf_k(l)$ cannot be
extendable to the cube $K^2$.

\smallskip This example shows that in the exceptional case $n=2, \ind =2$
the roots of holomorphic line bundles {\slsf do not extend} across
pseudoconcave boundaries. Indeed the critical case of index two can
be easily reduced to the model one, just discussed, via obvious
deformations and applying Theorem \ref{root-ext-thm}, see Section
\ref{sect.EXT} and Appendix for more details. It turns out and it is
our second result, that $\calf_k(l)$ exhibit {\slsf all}
nonextendable $k$-th roots.

\begin{nnthm}
\label{tors-ext-thm} In the conditions of Theorem \ref{root-ext-thm}
suppose that $n=2$ and $\ind\rho|_{\c}=2$. If $\calf$ doesn't extend
to a neighborhood of $\c$ then there exists a neighborhood $V$ of
$c$ biholomorphic to $K^2$ such that $\calf|_{K^2\cap U_{t_0}^+}$
extends to $D$ and is isomorphic to some  $\calf_k(l)$ there.
\end{nnthm}
This $\calf_k(l)$ is uniquely determined by $\calf$, see Remark
\ref{inv-l}.

\newprg[prgINT.contr]{Extension of roots across contractible analytic sets.}
 Now let us formulate a general result on
extending roots across contractible analytic sets. Recall that a
compact analytic set $E$ in a normal complex space $X$ is called
{\slsf contractible} if there exists a normal complex space $Y$, a
compact analytic set  $A$ in $Y$ of codimension at least two and a
holomorphic map (a contraction) $\C : X \to Y$ which is a
biholomorphism between $X\setminus E$ and $Y\setminus A$. If $A$ is
zero dimensional one calls $E$ {\slsf exceptional}.

\begin{nnthm}
\label{exc-ext-thm} Let $E$ be a contractible analytic set in a
normal complex space $X$ and let $\caln$ be a holomorphic line
bundle on $X$. Suppose $\caln$ admits a $k$-th root $\calf$ on
$X\setminus E$. Then $\calf$ extends to a holomorphic line bundle
$\widetilde{\calf}$ on the whole of $X$. Moreover,
$\widetilde{\calf}^{\otimes k}$ doesn't depend on $k$.
\end{nnthm}

\begin{nnrema} \rm
{\bf (a)} The extension $\widetilde{\calf}$ of $\calf$ can fail  to
be a $k$-th root of $\caln$ on $E$. Let $X$ be the blown-up $\cc^2$
in the origin  and $E$ be the exceptional curve. Set $\caln
\deff [E]$. By the adjunction formula $\caln|_E $ is the normal
bundle of the imbedding $E\subset X$ \ie is $\calo_E(-1)$. At the
same $\caln|_{X \setminus E}$ is trivial and as such admits a
trivial root $\calf =\calo$ on $X\setminus E = \cc^2\setminus
\{0\}$, and this $\calf$ is a root of $\caln|_{X\setminus E}$ of any
given degree $k\in \nn$.  But $\caln$ doesn't have roots of any
degree $k>1$ on $E$.

\smallskip\noindent{\bf (b)} At the same time Theorem \ref{exc-ext-thm}
means that $\caln|_{X\setminus E}$ can be (differently)  extended to
a line bundle, say $\widetilde{\caln}$ on $X$ having the extension
$\widetilde{\calf}$ of $\calf$ as its $k$-th root. And, moreover,
$\widetilde{\caln}$ is the same for all $k$. I.e., we can once
forever modify $\caln$ on the contractible set and make all roots of
$\caln$, which one can take on $X\setminus E$, to extend onto the
whole of $X$ as roots of the modified bundle.

\smallskip\noindent{\bf (c)} Note that according to our definition of
a contractible set all compact analytic sets of codimension at least
two are contractible. For them, in fact, more is true, see Lemma
\ref{space-ext}: $\calf$ extends to a $k$-th root of $\caln$ (and
not of some other $\widetilde{\caln}$).
\end{nnrema}

\begin{nnrema} \rm {\bf (a)} Among others extension results we prove
a Thullen type extension theorem for roots, see Theorem
\ref{thullen}.

\smallskip\noindent{\bf (b)} In section \ref{LOG} we remark that
analogous results can be obtained for {\slsf logarithms} of
holomorphic line bundles.

\smallskip\noindent{\bf (c)} We end up with an example of non-extendability
of roots from pseudoconcave in the sense of Andreotti domains in
complex projective plane.
\end{nnrema}

\newprg[prgINT.ack]{Acknowledgements.}{\bf 1.}
This paper in its final part was accomplished during the authors
stay in Math. Department of Nagoya University in September-October
2010. My stay was supported by Grant-in Aid for Scientific Research
(A) 22244008 and I'm grateful to Kengo Hirachi and Takeo Ohsawa both
for their hospitality and useful discussions. My interest in the
questions about extension of roots of holomorphic line bundles was
inspired by \cite{Oh}.

\smallskip\noindent{\bf 2.} I'm especially grateful to Stefan
Nemirovski for explaining to me an example of a non-extendable root
from a pseudoconvex in the sense of Andreotti domain in $\pp^2$.
This example is given in section \ref{LOG}.

\newsect[sect.EXT]{Extension of roots of holomorphic line bundles.}

\newprg[prgEXT.an-obj]{Extension of analytic objects} Throughout this
paper we shall use extension properties of various kind of "analytic
objects". These analytic objects have two decisive properties, which
we shall formalize in the  definition of  a {\slsf sheaf of analytic
objects} below.

\smallskip Let $\cala_X$ be a sheaf of sets on a complex  analytic
space $X$.

\begin{defi}
\label{sf-an-obj} A sheaf $\cala_X$ will be called a {\slsf sheaf of
analytic objects} if it satisfies the following two conditions:

\medskip{\slsf P1)} Sections of $\cala_X$ obey the uniqueness theorem.
I.e., if two of them $\sigma_1$, $\sigma_2$ are defined on a locally
irreducible, connected open subset $U\subset X$ and for some
nonempty open $V\subset U$ one has $\sigma_1|_V = \sigma_2|_V$, then
$\sigma_1=\sigma_2$ on $U$.

\medskip{\slsf P2)} The Hartogs type extension lemma is valid for
sections of $\cala_X$. I.e., if a section $\sigma$ of $\cala_X $ is
given over a Hartogs figure $E^n_{\eps} \deff H^2_{\eps}\times
\Delta^{n-2}\subset X$ then it extends to the section of $\cala_X$
over the polydisc $\Delta^n$, $n\ge 2$.

\medskip Sections of $\cala_X$ will be called {\slsf analytic objects}
over $X$.
\end{defi}
More precisely in (P2) we mean that this polydisc is imbedded as a
open subset $i(\Delta^n)$ into $X$ and $\sigma$ is initially defined
on $i(E^n_{\eps})$.  Recall that the two dimensional Hartogs figure
$H^2_{\eps}$ in this Definition is

\begin{equation}
\eqqno(hart-f1)
H^2_{\eps } = \Big(\Delta_{\eps}\times \Delta\Big)\bigcup \Big(\Delta
\times A_{1-\eps, 1}\Big),
\end{equation}
where $\Delta_{\eps}$ is the disc in $\cc$ of radius $\eps >0$,
$A_{1-\eps, 1} \deff \Delta\setminus \overline{\Delta}_{1-\eps}$ is
an annulus. Number $n\ge 2$ is the dimension of $X$ in a
neighborhood of $i(E^n_{\eps})$.

\begin{rema} \rm
\medskip\noindent{\bf (a)} The sheaf of holomorphic functions $\calo_X$ is a
sheaf of analytic objects.  The same with meromorphic functions:
$\calm_X = {\cal Mer} (X,\pp^1)$. Due to \cite{Iv2} instead of
$\pp^1$ one can take any compact K\"ahler $Y$.

\medskip\noindent{\bf (b)} A smooth holomorphic foliation on $X$ (resp.
singular holomorphic foliation) is defined by a holomorphic section
(resp. meromorphic section) of the appropriate Grassmann bundle $Gr
(TX)$. The products $U\times Gr (T_xX)\equiv Gr (TX)|_U$, where $U$
is a local chart, are K\"ahler (even projective). Therefore the
extension works, \ie the condition (P2) is satisfied.
Involutibility, being an analytic condition, is preserved by
extension. Uniqueness condition is also obvious. The extended
foliation might become singularities. Therefore the sheaf
$\calf_X(d)$ of codimension $d$ singular holomorphic foliations on
pure dimensional analytic space $X$ is a sheaf of analytic objects.

\medskip\noindent{\bf (c)} If a sheaf of objects satisfies the
property (P1) then the method of Cartan-Thullen provides the
existence of ''envelopes``, \ie of maximal domains over a the base
space $X$, where a given family of analytic objects extends. In fact
(P1) is not only sufficient but  also a necessary condition for
that, see \cite{Iv3}.

\medskip\noindent{\bf (d)} Coherent analytic sheaves/holomorphic
bundles are not analytic objects, they do not satisfy neither of
conditions (P1) and (P2).  This explains the problems with extending
them. The same with analytic sets. At the same time analytic
subsheaves of coherent analytic sheaves are analytic objects. This
was for the first time observed by Siu, see \cite{Si}, and morphisms
of bundles and sheaves are analytic objects (obvious).

\medskip\noindent{\bf (e)} Roots of holomorphic line bundles do
satisfy property (P2), see Lemma \ref{hart-l} below, but do not satisfy
(P1), see example in the Introduction.

\medskip\noindent{\bf (f)} Solutions of a general (\ie non-holomorphic)
elliptic system do satisfy (P1) but not (P2).
\end{rema}

We are going to state a general Bochner-Hartogs type extension
theorem for analytic objects. Recall that a smooth real valued
function $\rho$ on a complex manifold $U$ is called {\slsf strictly
plurisubharmonic} at point $x$ if the Levi form $dd^c\rho (x)$ is
positive definite. Here, as well as everywhere, by $U^+_t$ we denote
the superlevel set $\{ \rho (x) >t \}$ of $\rho$, by $\Sigma_t =
\{\rho =t\}$ the level set and by $U^-_t$ the lower level set
$\{\rho (x) <t\}$.

\begin{thm}
\label{ext-an-ob-thm} Let $W$ and $U$ be domains in a complex
manifold $X$ of dimension $n\ge 2$ and let $\rho$ be a strictly
plurisubharmonic Morse function on $U$ with $\rho (U) \supset
[t_1,t_0]$. Suppose that for every $t\in [t_1,t_0]$

\smallskip\sli the level set $\Sigma_t$ is connected,

\smallskip\slii the difference $\Sigma_t\setminus  W$ is compact,

\smallskip\sliii and that $U^+_{t_0} \subset W$.

\smallskip\noindent Then every analytic object extends from $W$ to
$U^+_{t_1}$.
\end{thm}
One can think  about $W$ as being a complement to a compact in $U$,
for example. The proof of this theorem is standard  and goes as
follows. Let $t$ be such that our analytic object $\sigma$ is
extended to $U_t^+$. All what one needs to do is to place Hartogs
figure appropriately near the hypersurface $\Sigma_t$ and, using
(P2) extend $\sigma$ through them, see Figure \ref{u1u2-fig}  for
the regular level sets. In the case of critical points one should
appropriately put Hartogs figures near the critical point. In
Appendix we give a precise Lemma \ref{index-2-lem}  of that kind.
The uniqueness condition (P1) will guarantee that all extensions
match one with another. And it is at this point where an analogous
Bochner-Hartogs type extension theorem fails for the roots of
holomorphic line bundles in the critical case - this will be
examined later on.

\newprg[prgEXT.gen]{Generalities on extensions of bundles and sheaves}
We refer to the book \cite{GR1} for the generalities on coherent
analytic sheaves and to \cite{Si} and \cite{ST} for the specific
questions about extensions of sheaves. Here we only briefly give a
few, very  basic, remarks about extensions of bundles and sheaves.
Bundles and their sheaves of sections will be always denoted with
the calligraphic letters, like $\caln$ or $\calf$. Their total
spaces - with the standard letters, like $N$ and $F$. If a line
bundle $\caln$ is canonically associated to a divisor $D$ we write
$\caln = \calo (D)$ and $N=[D]$ correspondingly.

\smallskip\noindent
{\bf 1.} Let $Y\subset \widetilde{Y}$ be domains in a complex
manifold (or, a normal complex space) $X$ and let $\caln$ be a
holomorphic bundle (resp. a coherent analytic sheaf) on $Y$. One
says that $\caln$ extends from $Y$ to $\widetilde{Y}$ if there
exists a holomorphic bundle $\widetilde{\caln}$ (resp. a coherent
analytic sheaf) on $\widetilde{Y}$ and an isomorphism $\phi :
\widetilde{\caln}|_Y\to \caln$ of bundles (resp. of sheaves).

\smallskip\noindent
{\bf 2.} If $\caln$ is a line bundle then it extends as a bundle if
and only if it extends as a coherent analytic sheaf. The way to see
this is to pass to the second dual $\widetilde{\caln}^{**}$ of the
extended sheaf, which is reflexive.  And then apply Lemma 26 from
\cite{Fr} to conclude that $\widetilde{\caln}^{**}$ is locally free,
\ie is a bundle extension of $\caln$.

\smallskip\noindent
{\bf 3.} The feature mentioned in the previous item is specific
for line bundles. The rank two subbundle
$\caln\subset \calo^3$ over $\cc^3\setminus\{0\}$ with the stalk
$\caln_z=\{ w: w_1z_1+w_2z_2+w_3z_3=0\}$ extends to the origin as a
coherent sheaf but not as a bundle.

\smallskip\noindent
{\bf 4.} Bundles/sheaves do not extend as
holomorphic/mero\-morphic functions, \ie they are not analytic
objects. The following example is very instructive. Take the
punctured bidisc $\Delta^2_* = \{z\in\cc^2:
0<\max{\{|z_1|,|z_2|\}}<1\}$ and cover it by two Stein domains
$U_1=\{z\in \Delta^2: z_1\not= 0\}$, $U_2=\{z\in \Delta^2: z_2\not=
0\}$. Every function of the form $f(z) =
e^{\sum_{k,l>0}\frac{a_{kl}}{z_1^kz_2^l}}$ is a transition function
of a holomorphic line bundle on $\Delta^2_*$, which doesn't extend
to the origin, provided at least one coefficient $a_{kl}\not=0$.
Indeed, it is easy to see that such a bundle doesn't admit a
non identically zero holomorphic section in a punctured neighborhood of the
origin. Moreover, all these bundles are pairwise non isomorphic.

In particular, we see that a Hartogs type extension condition (P2)
fails for holomorphic bundles. Uniqueness condition (P1) obviously
fails to, because all bundles of the same rank are locally
isomorphic to the trivial one.

\smallskip\noindent
{\bf 5.} At the same time  analytic singularities of codimension
three are removable for coherent sheaves. The rough reason is that
$H^{0,1}(\bb^3_*)=0$, see \cite{Si}. This last remark explains the
essential difference between dimension two and dimensions starting
from three in the problems addressed in this paper.

\smallskip\noindent
{\bf 6.}  Let us made the following:

\begin{rema} \rm
\label{ext-rem} {\bf (a)} Let a holomorphic  line bundle $\calf$ be
defined on a domain $V$ in a complex manifold $X$ and let $U$ be
another domain such that $U\cap V$ is nonempty and connected. If
$\calf|_{U\cap V}$ is trivial then $\calf$ extends to $V\cup U$.
Indeed, one can use the trivialization, say $\phi : \calf|_{U\cap
V}\to \cc\times (U\cap V)$ as a transition function to glue $\calf$
with the trivial bundle $\cc\times U$ on $U$.

\smallskip\noindent{\bf (b)} At the same time if $\calf$ is trivializable
on a subdomain $U_1\subset U\cap V$ then it might be not sufficient
to make such a gluing. We shall discuss such situation in the
Appendix.

\smallskip\noindent{\bf (c)} If there are two such domains $U_1$ and
$U_2$ and, in addition, $U_1\cap U_2$ and $U_1\cap U_2\cap V$ are
connected and if $\calf^i$ denotes the extensions  of $\calf$ onto
$U_i\cup V$, then we get a transition function $\phi_{12}: U_1\cap
U_2\cap V\to  \cc^*$. In order that $\calf^i$ glue together to a
holomorphic line bundle on $U_1\cup U_2\cup V$ it is necessary and
sufficient that $\phi_{12}$ extends to a non-vanishing function on
$U_1\cap U_2$. This last condition is not automatic.
\end{rema}

\begin{figure}[h]
\centering
\includegraphics[width=2.0in]{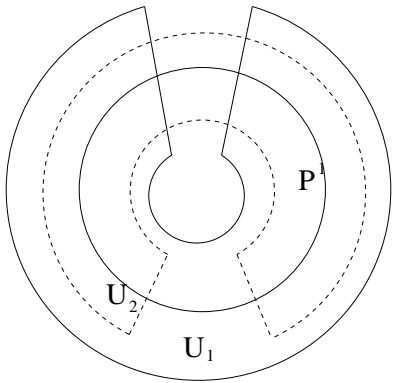}
\caption{The transition function $\phi_{12}$ between $\calf^{1}$ and
$\calf^{2}$ is defined only in $(U_1\cap U_2)\setminus \pp^1$. But it
should be defined on $U_1\cap U_2$ in order for $\calf$ to extend. }
\label{blow-fig}
\end{figure}

\begin{exmp} \rm
Consider the line bundle $\calf $ in $\Delta^2_*$ given by the
transition function $e^{\frac{1}{z_1z_2}}$. Blow up the origin and
denote by $\pp^1$ the exceptional divisor. Cover it by two bidiscs
$U_1$ and $U_2$ as on the Figure 1. Then $\calf|_{U_j\setminus
\pp^1}$ is trivial for $j=1,2$ (every holomorphic bundle on
$\Delta\times\Delta^*$ is trivial) and therefore extends to
$\Delta^2_*\cup U_j$ as a holomorphic line bundle $\calf^j$. But the
transition function between $\calf^1$ and $\calf^2$ doesn't extend
though $\pp^1\cap U_1\cap U_2$ and therefore $\calf$ doesn't extend
onto the blown up bidisc. And it shouldn't, because otherwise its
direct image would extend $\calf$ downstairs, and this is not the
case.
\end{exmp}

\newprg[EXT.thullen]{Hartogs Lemma and Thullen  type Extension Theorem for roots.}

Let us make the first step in the proof of Theorem
\ref{root-ext-thm} - a Hartogs type extension lemma for roots of
holomorphic line bundles. It occurs to be valid in all dimensions.
Recall that the Hartogs figure in $\cc^n$ is the following domain
\begin{equation}
\eqqno(hart-f2)
H^n_{\eps ,r} = \Big(\Delta^{n-1}_{\eps}\times
\Delta\Big)\bigcup \Big(\Delta^{n-1}\times A_{r , 1}\Big).
\end{equation}
Here  $\Delta\deff\Delta_1$, and $A_{r , 1} \deff
\Delta\setminus\bar\Delta_{r}$. When one takes $r=1-\eps$ one gets
the Hartogs figure from the Introduction. But in general the
possibility $r=0$ is not excluded (in that case $A_{0,1} = \Delta^*$
is the punctured disc), while $\eps$ should be strictly positive.

\newlemma{hart-l}
Let $\calf$ be a holomorphic line bundle on the Hartogs figure
$H^n_{\eps ,r}$, $n\ge 2$, and let $\caln=\calf^{\otimes k}$ be the tensor
$k$-th power of $\calf$, $k\ge 2$. Suppose that $\caln$ extends as a
holomorphic line bundle onto the unit polydisc $\Delta^n$. Then
$\calf$ extends as a holomorphic line bundle onto $\Delta^n$ and
stays there to be a $k$-th root of $\caln$.
\end{lem}
\proof Take the following covering  of $\Delta^n$: $U_1=
\Delta^{n-1}_{\eps}\times\Delta$ and $U_2=\Delta^{n-1}\times A_{r ,
1}$. Every holomorphic bundle on $U_j, j=1,2$ is trivial and
therefore $\calf$ is defined in $H^n_{\eps ,r}$ by a transition function
$f_{12}\in\calo^*(U_{12})$, here $U_{12}\deff U_1\cap U_2 =
\Delta^{n-1}_{\eps}\times A_{r , 1}$. Then $f_{12}^k$ is a
transition function of $\caln$. Every holomorphic bundle on
$\Delta^n$ is trivial and therefore there exist $F_1\in
\calo^*(U_1)$ and $F_2\in\calo^*(U_2)$ such that $f_{12}^k =
F_1\cdot F_2$ on $U_{12}$. $U_1$ is simply connected and therefore
$F_1$ admits a $k$-th root $f_1$, \ie such a non-vanishing
holomorphic function on $U_1$ that $f_1^k=F_1$. We have now that
$\Big(f_{12}\cdot f_1^{-1}\big)^k = F_2|_{U_{12}}$,\ie $F_2$ admits
a $k$-th root on $U_{12}=\Delta_{\eps}^{n-1}\times A_{r,1}$.
But then it obviously admits a $k$-th root on $\Delta^{n-1}\times
A_{r , 1} = U_2$. I.e., there exist $f_2\in\calo^*(U_2)$ such that
$\Big(f_{12}\cdot f_1^{-1}\big)^k = f_2^k$.  All what is left is to
multiply $f_2$ by an appropriate root of $1$ to get $f_{12}\cdot
f_1^{-1} = f_2$. Therefore $\calf$ is trivial and therefore extends to
$\Delta^n$.

\smallskip Let $\phi:\caln|_{H^n_{\eps ,r}} \to \calf^{\otimes k
}|_{H^n_{\eps ,r}}$ be an isomorphism.  In global holomorphic frames
$l$ of $\calf$ and $n$ of $\caln$ (over $\Delta^n$) $\phi$ is given
as $\phi(n) = h\cdot l^{\otimes k}$, where $h$ is a nonvanishing
holomorphic function in $H^n_{\eps ,r}$. By the classical Hartogs
extension theorem $h$ holomorphically extends to $\Delta^n$ and
doesn't vanish there. This gives the desired extension of $\phi$ and
proves that extended $\calf$ stays to be a $k$-th root of $\caln$.

\medskip\qed

\begin{rema} \rm
\label{non-uniq}
We see that the condition (P2) in the Definition \ref{sf-an-obj} of
analytic objects is satisfied by roots of holomorphic line bundles.
But it is easy to see that (P1) fails. Let $X$ be an Enriques
surface, \ie $\calk_X\not= \calo_X$ but $\calk_X^{\otimes 2} \cong
\calo_X$. Then for an imbedded bidisc $V\subset X$ we have that
$\calk_X|_V\cong \calo|_V$, but they are different globally. At the
same time both these bundles are square roots of the trivial bundle
over $X$.

\smallskip This was a {\slsf global} example. A local one was given in
the Introduction.
\end{rema}

When $r=0$ we obtain from Lemma \ref{hart-l} a Thullen type
extension lemma for the roots of holomorphic line bundles:

\begin{nncorol}
\label{thul-cor} Let $\calf$ be a holomorphic line bundle on
$(\Delta^{n-1}\times \Delta^*)\cup (\Delta^{n-1}_{\eps}\times
\Delta)$ and let $\caln=\calf^{\otimes k}$ be the tensor $k$-th
power of $\calf$, $k\ge 2$. Suppose that $\caln$ extends as a
holomorphic line bundle onto the unit polydisc $\Delta^n$. Then
$\calf$ extends as a holomorphic line bundle onto $\Delta^n$ and
stays there to be a $k$-th root of $\caln$.
\end{nncorol}

Let us remark that a more general Thullen type extension theorem for
roots of holomorphic line bundles is valid in all dimensions (\ie also in
dimension two).

\begin{thm}
\label{thullen} Let $A$ be an analytic set in a connected complex
manifold $X$ and let $G$ be a domain in $X$, which contains
$X\setminus A$ and which intersects every irreducible branch of $A$
of codimension one. Let $\caln$ be a holomorphic line bundle on  $X$
such that it admits a $k$-th root $\calf$ on $G$. Then $\calf$
extends to a holomorphic line bundle onto the whole of $X$ and stays
there to be a $k$-th root of $\caln$.
\end{thm}
\proof Let us establish the uniqueness of the extension first.
Suppose $\calf$ is extended as a $k$-th root of $\caln$ to a bundle
$\calf_1$ on an open set $G_1\supset G$ and, again as a $k$-th root
of $\caln$ to a bundle $\calf_2$ on $G_2\supset G$. Fix an
isomorphism $H:\calf_1|_G\to \calf_2|_G$. Let us see that $H$
extends to an isomorphism between $\calf_1|_{G_1\cap G_2}$ and
$\calf_2|_{G_1\cap G_2}$ and therefore $\calf$ extends to $G_1 \cup
G_2$ as a $k$-th root of $\caln$.

\smallskip Denote by $\Reg A$ the regular part of the codimension one
locus of $A$. Fix a point $q\in \Reg A\cap G$ and take an another
point $p\in \Reg A\cap (G_1\cap G_2)$ on the same irreducible
component of $\Reg A$ as $q$. Take now open subsets $q,p\in
U_1\subset \Reg A$ and $U\subset X$ such that

\smallskip\sli $U_1=\Delta^{n-1}$ with $q=0$ and $p=(1/2,0,...,0)$;

\smallskip\slii $U=\Delta^n$ and $U\cap \Reg A = U_1$;

\smallskip\sliii $\Delta^{n-1}_{\eps}(q)\times \Delta\subset G$ and
$\Delta^{n-1}_{\eps}(p) \times \Delta\subset G_1\cap G_2$ for some
$\eps >0$.

\smallskip Such $U_1, U$ can be found using the Royden's lemma, see \cite{Ro}.
Set $V \deff (U\setminus U_1)\cup (\Delta^{n-1}_{\eps}(q)\times \Delta) \cup
(\Delta^{n-1}_{\eps}(p)\times \Delta )$.

\smallskip\noindent{\slsf Claim. $\calf_i$ admits a global holomorphic frame
over $V$.} Indeed, set $V_1 \deff (U\setminus U_1)\cup
(\Delta^{n-1}_{\eps}(q) \times \Delta )$ and $V_2 \deff (U\setminus
U_1)\cup (\Delta^{n-1}_{\eps}(p) \times \Delta )$. By Lemma
\ref{hart-l} $\calf_i$ admits global holomorphic frames $g_j$ over
$V_j$, $j=1,2$. We have that $g_1 = hg_2$ on $U\setminus U_1$ for
some nonvanishing holomorphic function $h$ on $U\setminus U_1$. Let
$\psi_i:\calf_i^{\otimes k}\to \caln|_{G_i}$ be an isomorphism.
Remark that $\psi_i(g_j^{\otimes k})$ for $j=1,2$ are global holomorphic frames
of $\caln|_{V_j}$ and by the standard Hartogs theorem they extend
holomorphically onto $U$ as frames of $\caln$. From this fact we see
that $h^k$ and therefore $h$ extend to a nonvanishing holomorphic
function on $U$. I.e., $hg_2$ is a global holomorphic frame of
$\calf_i$ over $V$. The Claim is proved.

\smallskip Let $f_i$ be global frames of $\calf_i$ over $V$ for $i=1,2$,
constructed in the Claim above. Since $V_1\subset G$ we see that
there exists a nonvanishing holomorphic function $h$ on $V_1$ such that
$H(f_1) = hf_2$. Again by Hartogs theorem $h$ extends to a nonvanishing
holomorphic function on $U$. Therefore we can extend $H$ onto $V_2$
by setting  $H(f_1|_{V_2}) = hf_2|_{V_2}$.

\smallskip We extended $H$ through every regular point of $A$. Extension of
$H$ through the singular locus of $A\cap G_1\cap G_2$ goes along the
same lines using the stratification of $\Sing A$. Uniqueness of the
extension is proved.

\smallskip Having established the uniqueness of the extension we can finish
the proof of our Theorem with the help of Corollary \ref{thul-cor}.
Let $\tilde G\supset G$ be the maximal domain onto which $\calf$
extends as a $k$-th root of $\caln$. If there exists $p \in \d\tilde
G\cap \Reg A$ then we can apply Corollary \ref{thul-cor} (with
$r=0$) and extend $\calf$ to a neighborhood of $p$. This contradicts
to the maximality of $\tilde G$. The removal of $\Sing A$ is
analogous.

\smallskip The condition on the extended bundle to be a $k$-th root of $\caln$
is also preserved by the Thullen type theorem for holomorphic functions.

\smallskip\qed

\begin{rema} \rm
We place here a Thullen type extension theorem for roots  for the
following reasons:

\smallskip\noindent{\bf (a)} In the situation when roots of line bundles
are not always extendable it is desirable to have complete
understanding when they nevertheless  do extend.

\smallskip\noindent{\bf (b)} In the process of the proof it was very
clear how the uniqueness property plays a crucial role in extension.

\smallskip\noindent{\bf (c)} The proof of the extendability through
subcritical points in dimension two and through critical one in
dimension $\ge 3$ is analogous to that given for Thullen case.
\end{rema}

\newprg[EXT.reg]{Extension of roots across regular level sets.}
We shall make a first step in the proof of Bochner-Hartogs type
Extensions Theorems \ref{root-ext-thm} and \ref{tors-ext-thm}  for
roots of holomorphic line bundles. For the convenience of the future
references we shall state Theorems \ref{root-ext-thm} and
\ref{tors-ext-thm} in a slightly more general, but obviously
equivalent form. As in Theorem \ref{ext-an-ob-thm} we consider two
domains $W$ and $U$ in a complex manifold $X$ of dimension $n\ge 2$
and let $\rho$ be a strictly plurisubharmonic Morse function on $U$
with $\rho (U) \supset [t_1,t_0]$. We shall suppose that for all
$t\in [t_1,t_0]$ the level set $\Sigma_t$ is connected,
$\Sigma_t\setminus W$ is compact and that $U_{t_0}^{+}\subset W$.

\begin{thm}
\label{gen-ext-thm} Let $\caln$ be a holomorphic line bundle on $X$
and suppose that it admits a $k$-th root $\calf$ on $W$ for some
$k\ge 2$. Then:

\smallskip\sli $\calf$ extends to a $k$-th root of $\caln$ on $U_{t_1}^{+}$
provided $n=\dim U\ge 3$ or, $n=2$ and neither of $\Sigma_t$-s has
critical points of index $2$.

\smallskip\slii If $n=2$, $\Sigma_{t^*}$ for some $t^*<t_0$ contains some
critical point $\c$ of index $2$, and $\calf$ is extended to
$U_{t^*}^{+}$ as a $k$-th root of $\caln$ then either $\calf$
extends to a neighborhood of $\c$ or, there exists a neighborhood
$V$ of $\c$ biholomorphic to $K^2$ such that $\calf|_{K^2\cap
U_{t^*}^+}$ extends to $D$ and is isomorphic there to some uniquely
defined $\calf_k(l)$.
\end{thm}
\proof By the assumption $\calf$ exists as a $k$-th root of
$\caln|_{U_{t_0}^{+}}$ on $U_{t_0}^{+}$. Let $\psi_{t_0}:
\caln|_{U_{t_0}^{+}} \to \calf^{\otimes k}$ be a corresponding
isomorphism. Denote by $t^*$ the infinum of such $t\le t_0$  that
$\calf$ extends up to $t$ as a $k$-th root $\calf_t$ of
$\caln|_{U_t^{+}}$. More precisely, by saying that $\calf$ extends
up to $t$ as a $k$-th root $\calf_t$ of $\caln$, we mean the
following:

\medskip a) Holomorphic line bundles $\calf_{t_1}$ on $U_{t_1}^{+}$ are
defined for all $t_1 \ge t$ together with isomorphisms
$\psi_{t_1}:\caln|_{U_{t_1}^{+}}\to \calf^{\otimes k}_{t_1}$.

\smallskip b) For all pairs $t_1 \ge t_2\ge t$ the bundle $\calf_{t_2}$
is an extension of $\calf_{t_1}$, \ie an isomorphism
$\phi_{t_1t_2}:\calf_{t_2}|_{U_{t_1}^{+}} \to \calf_{t_1}$ is given,
and these isomorphisms satisfy:

\smallskip\qquad $\text{b)}_1$ $\phi_{t_1t_1}=\id$;

\smallskip\qquad $\text{b)}_2$ $\phi_{t_1t_2}\circ \phi_{t_2t_3}=\phi_{t_1t_3}$
for $t_1\ge t_2 \ge t_3$.

\smallskip c) Isomorphisms $\psi_{t_1}$ and $\phi_{t_1t_2}$ are natural in
the sense that for every pair $t_1\ge t_2\ge t$ the following
diagram is commutative:
\begin{equation}
\def\normalbaselines{\baselineskip20pt\lineskip3pt \lineskiplimit3pt }
\def\mapright#1{\smash{\mathop{\longrightarrow}\limits^{#1}}}
\def\mapdown#1{\Big\downarrow\rlap{$\vcenter{\hbox{$\scriptstyle#1$}}$}}
\begin{array}{ccccccccc}
\caln|_{U_{t_{2}^{+}}}&\mapright{r_{t_1t_2}}&\caln|_{U_{t_{1}^{+}}}&\\
\mapdown{\psi_{t_{2}}}& &\mapdown{\psi_{t_{1}}}& \\
\calf^{\otimes k}_{t_{2}}&\mapright{\phi^{\otimes k}_{t_1t_{2}}}&
\calf^{\otimes k}_{t_{1}}& .\\
\end{array}
\eqqno(proj-0)
\end{equation}
Here $r_{t_1t_2}:\caln|_{U_{t_{2}^{+}}}\to \caln|_{U_{t_{1}^{+}}}$
is the natural restriction operator of the globally existing bundle $\caln$.

\medskip\noindent{\slsf Step 1. The infinum $t^*$  is achieved.}
If $t^*$ is the infinum then for every $t>t^*$ the bundle $\calf$
extends onto $U^{+}_t$ as a $k$-th root $\calf_t$ of $\caln$. We define the
presheaf $\calf_{t^*}$ on $U^{+}_{t^*}$ as the projective limit of
$\calf_t$-s for $t>t^*$:
\begin{equation}
\eqqno(proj-1)
\calf_{t^*} \deff \lim_{\longleftarrow }\calf_t.
\end{equation}
A section $\sigma$ of $\calf_{t^*}$ over an open set $V\subset
U^{+}_{t^*}$ is a product $\prod_{t>t^*}\sigma_t$ of sections $\sigma_t$
of $\calf_t$ over $V\cap U^{+}_t$ such that $\phi_{t_1t_2}(\sigma_{t_2})
 =\sigma_{t_1}$ for every pair $t_1\ge t_2>t^*$.
Restriction map for $W\subset V$ in $\calf_{t^*}$ is
$\prod_{t>t^*}\sigma_t \to \prod_{t>t^*}\sigma_t|_{W\cap U^{+}_t}$,
which is correctly defined.

\smallskip It is easy to see that the presheaf, so defined, is actually a
sheaf and that this sheaf is locally free of rank one.  Isomorphisms
$\phi_{tt^*}: \calf_{t^*}|_{X_t}\to \calf_t$ for $t>t^*$ are
naturally defined as $\phi_{tt^*} (\prod_{t>t^*}\sigma_t) =
\sigma_t$. For $t=t^*$ we set $\phi_{t^*t^*}=\id$.

\smallskip Take some $t_1>t^*$. It is not difficult now to see that $\psi_{t_1}:
\caln|_{U_{t_1}^{+}} \to \calf_{t_1}^{\otimes k}$ extends to an isomorphism
$\psi_{t^*}:\caln|_{U_{t^*}^{+}} \to \calf_{t^*}^{\otimes k}$ and the
following diagram is commutative:

\begin{equation}
\def\normalbaselines{\baselineskip20pt\lineskip3pt \lineskiplimit3pt }
\def\mapright#1{\smash{\mathop{\longrightarrow}\limits^{#1}}}
\def\mapdown#1{\Big\downarrow\rlap{$\vcenter{\hbox{$\scriptstyle#1$}}$}}
\begin{array}{ccccccccc}
\caln|_{U_{t^*}^{+}}&\mapright{r_{t_1t^*}}&\caln|_{U_{t_1}^{+}}&\\
\mapdown{\psi_{t^*}}& &\mapdown{\psi_{t_1}}&\\
\calf^{\otimes k}_{t^*}&\mapright{\phi^{\otimes k}_{t_1t^*}}&
\calf^{\otimes k}_{t_1}. \\
\end{array}
\eqqno(proj-2)
\end{equation}
The proof goes by continuity: for, if $\psi_{t_1}$ is
already extended to $\psi_{t}$, then take the composition
$(\phi_{tt^*}^{-1})^{\otimes k}\circ \psi_{t}: \caln|_{X_{t}} \to
\calf_{t^*}^{\otimes k}|_{X_{t}}$. It is an analytic morphism of sheaves
and in local chats it is given by holomorphic functions. Therefore
it can be extended through a pseudoconcave boundary.

\begin{rema} \rm
\label{non-reg}
On this step we do not need to suppose that $\Sigma_{t^*}$ is regular.
\end{rema}

 In the remaining steps we shall bypass $t^*$. In the process
of the proof we shall distinguish the case when $\Sigma_{t^*}$ is  regular
 from the case when it is a critical level of some strictly plurisubharmonic
 Morse function, say $\rho$.

\smallskip\noindent{\slsf Step 2. Passing through a regular value.}
Let $t^*$ be a regular value of $\rho$ and suppose that
$\calf_{t^*}$ is an extension of $\calf$ to $U^{+}_{t^*}$, which stays
to be a $k$-th root of $\caln$ there. Take a point $z_0\in
\Sigma_{t^*}$. Then in an appropriate coordinates near $z_0$ one can
find a polydisc $\Delta^n$ such that $\Delta^n\cap U^{+}_{t^*}$ is
connected, $z_0=(\frac{1}{2},0,...,0)$ and such that
$H^n_{\eps}\subset U^{+}_{t^*}$ for some $\eps >0$. By Lemma
\ref{hart-l} the restriction $\calf_{t^*}|_{H^n_{\eps}}$ is trivial;
let $\phi : \calf_{t^*}|_{H^n_{\eps}}\to \cc\times H^n_{\eps}$ be a
trivialization. Extend it by the usual Hartogs theorem to an
isomorphism $\phi : \calf_{t^*}|_{\Delta^n\cap U^{+}_{t^*}}\to
\cc\times\Big(\Delta^n\cap U^{+}_{t^*}\Big)$. Then it defines (as a
transition function) a bundle $\tilde\calf$ over $U^{+}_{t^*}\cup
\Delta^n$ - an extension of $\calf_{t^*}$. Furthermore, by
assumption there exists an isomorphism $\psi : \calf^{\otimes
k}_{t^*}|_{U^{+}_{t^*}} \to \caln|_{U^{+}_{t^*}}$. The isomorphism
$\psi\circ (\phi^{-1})^{\otimes k}:\tilde \calf|_{\Delta^n\cap
U^{+}_{t^*}}^{\otimes k}\to \caln|_{\Delta^n\cap U^{+}_{t^*}}$ obviously
extends to $\Delta^n$ (again by the classical Hartogs theorem) and
this gives an extension of $\psi$ onto $U^{+}_{t^*}\cup \Delta^n$.
Now we can cover $\Sigma_{t^*}$ by a finite number of such
coordinate neighborhoods $U_j$ that:

\medskip
\sli Each $U_j$ is biholomorphic to a polydisc $\Delta^n$ and
$U_j\cap U^{+}_{t^*}$ are connected;

\smallskip
\slii the associated Hartogs figures $H^n_{\eps}$ are contained in
$U_j\cap U^{+}_{t^*}$;

\smallskip
\sliii the double $U_{ij} \deff U_i\cap U_j$ and triple
$U_{ijk}\deff U_i\cap U_j\cap U_k$ intersections are connected,
simply connected and moreover, $U_{ij}\cap U^{+}_{t^*}$ and $U_{ijk}\cap
U^{+}_{t^*}$ are connected and simply connected as well.

\medskip Extend $\calf_{t^*}$ through all of them using Lemma \ref{hart-l},
\ie get line bundles $\calf_{t^*}^j$ over $V_j\deff U^{+}_{t^*}\cup U_j$
- extensions of $\calf_{t^*}$. We need to prove that these
extensions match together, \ie that transition functions extend from
$U_i\cap U_j\cap U^{+}_{t^*}$ to $U_i\cap U_j$ for all $U_i\cap
U_j\not=0$. Denote by $\phi_j:\calf_{t^*}^j|_{U^{+}_{t^*}}\to
\calf_{t^*}$ the corresponding isomorphisms. Then we get the
transition functions $\phi_{ij}=\phi_i^{-1}\circ \phi_j:
\calf_{t^*}^j|_{U_i\cap U_j\cap U^{+}_{t^*}} \to \calf_{t^*}^i|_{U_i\cap
U_j\cap U^{+}_{t^*}}$. As it was said already, we need to extend
$\phi_{ij}$ onto $U_i\cap U_j$ (this is not automatic, see Fig. \ref{u1u2-fig}).

\begin{figure}[h]
\centering
\includegraphics[width=2.5in]{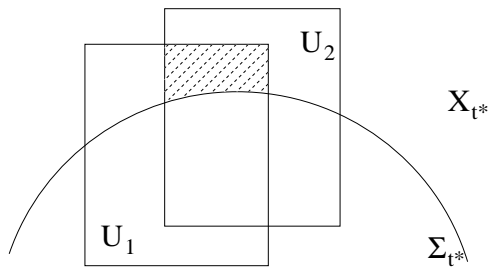}
\caption{The transition function $\phi_{12}$ between
$\calf_{t^*}^{1}$ and $\calf_{t^*}^{2}$ is defined only in $U_1\cap
U_2\cap U_{t^*}^{+}$ - dashed zone on the picture. But it should be
defined on $U_1\cap U_2$. Since the envelope of holomorphy of
$U_1\cap U_2\cap U^{+}_{t^*}$ is much smaller then $U_1\cap U_2$ the
extension of $\phi_{12}$ to $U_1\cap U_2$ is not automatic, \ie
cannot be achieved via the classical Hartogs theorem.}
\label{u1u2-fig}
\end{figure}

\smallskip Denote by $\psi_j:\caln|_{V_j}\to (\calf_{t^*}^j)^{\otimes k}$
the isomorphisms, obtained in Lemma \ref{hart-l} as extensions of
the isomorphism $\psi_{t^*}:\caln|_{U_{t^*}^{+}}\to \calf_{t^*}^{\otimes
k}$. From the diagram

\begin{equation}
\def\normalbaselines{\baselineskip20pt\lineskip3pt \lineskiplimit3pt }
\def\mapright#1{\smash{\mathop{\longrightarrow}\limits^{#1}}}
\def\mapdown#1{\Big\downarrow\rlap{$\vcenter{\hbox{$\scriptstyle#1$}}$}}
\begin{array}{ccccccccc}
\caln|_{U_i\cap U_j\cap U^{+}_{t^*}}&\mapright{\tr}&\caln|_{U_i\cap U_j\cap U^{+}_{t^*}}&\\
\mapdown{\psi_j}& &\mapdown{\psi_i}&\\
(\calf_{t^*}^j|_{U_i\cap U_j\cap U^{+}_{t^*}})^{\otimes
k}&\mapright{\phi^{\otimes k}_{ij}}&
(\calf_{t^*}^i|_{U_i\cap U_j\cap U^{+}_{t^*}})^{\otimes k} \\
\end{array}
\eqqno(match1)
\end{equation}
we see that $\phi^{\otimes k}_{ij} = \psi_i\circ\psi_j^{-1}$ on $U_i\cap U_j\cap U^{+}_{t^*}$. But
$\psi_j$ (resp. $\psi_i$) is defined over $U_j$ (resp. $U_i$) and $\tr$ is a transition map of a globally
existing bundle $\caln$, therefore $\psi_i\circ\tr\circ\psi_j^{-1}$ is defined over  $U_i\cap U_j$ extending
$\phi^{\otimes k}_{ij} $. But then $\phi_{ij}$ also extends onto $U_i\cap U_j$ as a $k$-th root on an
extendable non-vanishing function on a simply connected domain.

\smallskip
It is easy to see that the cocycle condition for the extended
transition maps will be preserved. Indeed, it is satisfied on
$U_{ijk}\cap U^{+}_{t^*}$, so it will be satisfied on $U_{ijk}$ to.

\newprg[EXT.chern]{Line bundles with vanishing Chern class.}
Before considering the case of critical points let us make few
preparatory remarks. Let $X$ be a connected (not necessarily
compact) complex manifold. Then we have the following exact
sequence:

\begin{equation}
\eqqno(jakobi1)
0\to \zz \to H^0(X, \calo) \buildrel \exp \over  \longrightarrow H^0(X, \calo^*) \to  H^1(X,\zz)
\to H^{0,1}(X) \to Pic^0(X)\to 0,
\end{equation}
where $Pic^0(X)$ is the group of topologically trivial holomorphic
line bundles on $X$. $Pic^0(X)$ appears here as the kernel of the
map $c_1$ from the group $H^1(X,\calo^*)$ of all  holomorphic line
bundles on $X$ to the N\'eron-Severi group $NS(X)$. I.e., $Pic^0(X)$
is exactly the group of holomorphic line bundles on $X$ with
vanishing first Chern class. The arrow  $H^{0,1}(X) \to Pic^0(X)$ is
the composition of the Dolbeault isomorphism $D: H^{0,1}(X)\to
H^1(X, \calo)$ and exponential map. If the map
\[
H^0(X,\calo)\buildrel \exp \over \to H^0(X,\calo^*)
\]
is surjective  then \eqqref(jakobi1)writes as

\begin{equation}
\eqqno(jakobi2)
0 \to  H^1(X,\zz) \to H^{0,1}(X) \to Pic^0(X)\to 0,
\end{equation}
which means  that  $Pic^0(X)\cong H^{0,1}(X)/H^1(X,\zz)$. This is
the case for example if all holomorphic functions on $X$ are
constant or, all holomorphic functions from $X$ extend to some
simply connected $\widetilde{X}\supset X$. The pair $K^2\supset
D\deff K^2\setminus \{x_1=x_2=0\}$ from the example in the
Introduction is the case we are particularly going to exploit.

\newprg[EXT.sub-crit]{Extension through subcritical points.}

Now we shall consider the case with critical points. Let us make
first of all some standard observations. We follow computations from
\cite{Iv1}, though other sources also might be used at this place.
Let $t^*$ be a critical value of our strictly plurisubharmonic Morse
function $\rho$ and let $z_0$ be a critical point of $\rho$ on
$\Sigma_{t^*} \deff \{ \rho(x) = t^*\}$. In an appropriate
coordinates near $z_0=0$ (and make $t^*=0$ for convenience)  the
function $\rho$  has the form

\begin{equation}
\eqqno(rho)
\rho (z) = \sum_{j=1}^na_jz_j^2 + \sum_{j=1}^na_j\bar z_j^2 + \sum_{j=1}^n
\vert z_j\vert^2 + O(\Vert z\Vert^3)
\end{equation}
with some real non-negative coefficients $a_1, ..., a_n$. In
coordinates $z_j=x_j+iy_j$ we rewrite \eqqref(rho) as follows

\begin{equation}
\eqqno(rho+)
\rho (z) = 2\sum_{j=1}^na_j(x_j^2-y_j^2) + \sum_{j=1}^n(x_j^2+y_j^2) +
O(\Vert z\Vert^3) =
\end{equation}
\[
= \sum_{j=0}^p\Big[(1+2a_j)x_j^2 + (1-2a_j)y_j^2\Big] +
\sum_{j=p+1}^n \Big[(1+2a_j)x_j^2+(1-2a_j) y_j^2\Big] + O(\Vert
z\Vert^3),
\]
where the index $p$ has the following significance:
$a_j>\frac{1}{2}$ for $j\le p$ and $0\le a_j<\frac{1}{2}$ for $j\ge
p+1$. Number $p$ is in fact the index of the Morse function $\rho$
at $0$. Case $a_j=\frac{1}{2}$ is excluded by the choice of $\rho$
to be Morse, case $p=0$ is not excluded. Taking $\delta >0$ and
$\norm{z}$ sufficiently small we can write
\begin{equation}
\eqqno(rho-1) \rho (z) \ge \sum_{j=0}^p[(2a_j+1)x_j^2-(2a_j-1)y_j^2]
+ \delta\cdot \sum_{j=p+1}^n\vert z_j\vert^2 + O(\Vert z\Vert^3) \ge
\end{equation}
\[
\ge a^2\sum_{j=0}^px_j^2 - b^2\sum_{j=0}^py_j^2 + \delta^2\cdot
\sum_{j=p+1}^n\vert z_j\vert^2  := \rho_1(z)
\]
with some $a>b$. From \eqqref(rho-1) we see that $Y^+\deff \{ z\in
\Delta^n: \rho_1(z)>0\} \subset U_{t^*}^{+}$ near the origin.
Therefore all that we need is to extend the root $\calf$ from $Y^+$
to $Y^+\cup \{\text{a neighborhood of the origin}\}$.  Denote the
coordinates $(z_1,...,z_p)$ as $z$ with $z=x+iy$ in vector
notations, and $(z_{p+1},...,z_n)$ by $w$. Set $A =\frac{a}{b}>1$,
$K_1=\{x\in \rr^p: \norm{x}< 1\}$, $K_A=\{y\in\rr^p: \norm{y}< A\}$
and $K^2=K_1\times K_A\times\Delta^{n-p}$. Appropriately deforming
$\Sigma_{t^*}^{+}$ and rescaling coordinates near the origin  we can
suppose that
\begin{equation}
\eqqno(y+1)
Y^+\deff \{(z,w)\in K^2: \rho_1(z) = a^2\norm{x}^2 - b^2\norm{y}^2
+ \delta^2\norm{w}^2>0\}.
\end{equation}

\smallskip\noindent{\slsf Sep 3. Subcritical case, \ie $0\le p\le n-1$.}
The case $p=0$ is obviously served by Hartogs type Lemma
\ref{hart-l}.

\begin{lem}
\label{sub-crit}
In the case $1\le p \le n-1$ domain $Y^+$ is simply connected.
\end{lem}
\proof  This is immediate. Indeed, from \eqqref(y+1) we see that for
$z$ such $a^2\norm{x}^2 - b^2\norm{y}^2>0$ the point $w$, with
$(z,w)\in Y^{+}$, belongs to an appropriate $(n-p)$-disc; if
$a^2\norm{x}^2 - b^2\norm{y}^2=0$ then $w$ belongs to a punctured
$(n-p)$-disc; and, finally, if  $a^2\norm{x}^2 - b^2\norm{y}^2<0$
then $w$ lies in an appropriate ring domain of the type
$A^{n-p}_{r_1,r_2} \deff
\Delta_{r_1}^{n-p}\setminus\overline{\Delta}_{r_2}^{n-p}$.  It is
important here that $n-p>0$. Now the $w$-coordinate of every loop in
$Y^+$ can be moved to this ring domain and then $z$-coordinate
retracted to a point. One then moves this point to a point inside
$\{a^2\norm{x}^2 - b^2\norm{y}^2<0\}$ and finally  contracts the
loop.

\smallskip\qed

\smallskip Since $\caln$ is extendable from $Y^+$ to $K^2$ we have that
$c_1(\caln|_{K^2})=0$ and therefore $c_1(\calf|_{Y^+})=0$. As the
result $\calf|_{Y^+}$ belongs to $H^{0,1}(Y^+)$. Let $\omega$ be a
$(0,1)$-from in $Y^+$ representing $\calf|_{Y^+}$ and let $\tau$ be
a $(0,1)$-from in $K^2$ representing $\caln|_{K^2}$. Then, since
$\calf^{\otimes k}$ extends to $K^2$ and is equal to $\caln$ there,
and because $H^1(Y^+,\zz)=0$ one has
\[
\tau|_{Y^+} = k\omega + \dbar f
\]
for some smooth function $f$ in $Y^+$. But $\caln|_{K^2}$ is trivial
and therefore $\tau = \dbar g$ for some smooth function $g$ on
$K^2$. Therefore
\[
\omega = 1/k(\dbar g - \dbar f),
\]
which means that $\calf $ is trivial on $Y^+$ to. Now one can extend
$\calf$ to $K^2\cup U_{t^*}^{+}$ as in Remark \ref{ext-rem}. The
case is proved.

\newprg[EXT.crit]{Critical case.}

Now we want to bypass critical points of maximal index $p=n$.  In
that case we have that  $K^2=K_1\times K_A$ is an open brick in
$\cc^n$ with $K_1=\{x\in \rr^n: \norm{x}< 1\}$, $K_A=\{y\in\rr^n:
\norm{y}< A\}$ and
\begin{equation}
\eqqno(y-plus-f)
Y^+\ = \bigl\{ z\in K^2: a^2\norm{x}^2>b^2\norm{y}^2 \bigr\} = \bigl\{z\in K^2:
A^2\norm{x}^2>\norm{y}^2\bigr\}.
\end{equation}
Consider one more strictly plurisubharmonic function
\begin{equation}
\eqqno(rho-2)
\rho_2(z) = A^2\norm{x}^2 -\norm{y}^2.
\end{equation}

\smallskip\noindent{\slsf Step 4. Critical case in dimension $\ge 3$.}
In that case $Y^+$ is again obviously simply connected and the same
proof as in subcritical case gives an extension of $\calf$ to a
neighborhood of the critical point.

\smallskip Theorem \ref{root-ext-thm} together with the part (\sli of
Theorem \ref{gen-ext-thm} are proved.

\medskip Now let us make a preparatory remark for the proof of Theorem \ref{tors-ext-thm}.
 Denote by $Y_{\tau}\deff \{z\in K^2:\rho_2 (z)>\tau \}$ - the
superlevel set of $\rho_2$ and by $S_{\tau}\deff \{z\in K: \rho_2
(z)=\tau \}$ denote the level set. Then $Y^+ = Y_0$ in our new
notations.
\begin{rema} \rm
Using the fact that $A>1$ one can deform $S_{\tau}$ to hypersurfaces
$\tilde S_{\tau}$ in such a way that they will stay smooth and
strictly pseudoconvex for all $\tau >0$ and moreover, they will
exhaust the following domain $D_{\delta,\eps}$ for an appropriate
$0<\eps<<\delta << 1$. Namely, set $K_{\delta} \deff \{(x_1,x_2):
-\delta < x_1, x_2 < \delta\}$ and  $K_{\eps} \deff \{(y_1,y_2):
-\eps < y_1, y_2 < \eps\}$.  Then set also  $K^2_{\delta,\eps}
\deff  K_{\delta}\times K_{\eps}$, and define $D_{\delta,\eps} \deff K^2_{\delta,\eps}\setminus \{x_1= x_2= 0\}$.
 \end{rema}

\newprg[EXT.twist-ext]{Case of index and dimension two}
Let us start this Subsection with a remark about torsion bundles.
Denote by $\zz_k = \{e^{\frac{2\pi il}{k}}: 0\le l\le k-1\}$
the group of $k$-th roots of $\adyn$.

\begin{rema} \rm
\label{tors-rem} Let $\cale$ be a torsion bundle of order $k$ on a
manifold $V$. Let us see that it can be defined by a locally
constant cocycle $\eps_{jl} \in \zz_k$ in an appropriate fine
covering $\{U_j\}$ of $V$. Indeed, let $\{U_j\}$ be a Stein covering
with simply connected $U_j$ and let $e_{jl}$ be a defining cocycle
for $\cale$ in this covering. Since $\cale^{\otimes k}$ is trivial
there exist holomorphic nonvanishing $g_j$ in $U_j$ such that
$g_je_{jl}^kg_l^{-1} = \adyn$. Take $k$-th roots $\tilde g_j$ of
$g_j$ on $U_j$ arbitrarily. Then $\tilde e_{jl} \deff \tilde g_j
e_{jl}\tilde g_l^{-1}$ will be still a defining cocycle for $\cale$.
But now $\tilde e_{jl}^k = \adyn$, \ie $\eps_{jl} \deff \tilde
e_{jl} \in \zz_k$.
\end{rema}

We are prepared to give the proof of Theorem \ref{tors-ext-thm}. Let
$\c \in \Sigma_{t_0}$ be a critical point in question. Choosing
coordinates around $\c=0$ appropriately and deforming $\Sigma_t$ for
$t>t_0$ (equivalently $Y_{\tau}$ for $\tau >0$) near $\c$, as it was
explained in the previous subsection, we can achieve, using Theorem
\ref{root-ext-thm} the extendability of $\calf$ onto
$D_{\delta,\eps}\deff K_{\delta}\times K_{\eps}\setminus \{x_1= x_2=
0\}$ in these coordinates.

\smallskip Cover $D_{\delta,\eps}$ by four Stein domains $D_1\deff \{x_1>0\}
\cap D_{\delta,\eps}$, $D_2 \deff \{x_2>0\}\cap D_{\delta,\eps}$,
$D_3\deff \{x_1<0\}\cap D_{\delta,\eps}$, $D_4 \deff \{x_2<0\} \cap
D_{\delta,\eps}$. For naturals $1\le l <k$ a holomorphic line bundle
$\calf_k(l)$ on $D_{\delta,\eps}$ from Introduction is defined by
the following cocycle $\{f_{ij}\} \in H^1(\calo^*, \{D_i\})$:

\begin{equation}
\eqqno(bund-f)
f_{ij} =
\begin{cases}
e^{\frac{2\pi il}{k}} & \text{ if } i=1, j=2;\cr 1 &\text{ otherwise
}.
\end{cases}
\end{equation}

\begin{lem}
\label{torsion} Let $\calf$ be a torsion bundle on $D_{\delta,\eps}$
of order $k\ge 2$. Then $\calf$ is isomorphic to one of
$\calf_k(l)$.
\end{lem}
\proof Since $H^2(D_{\delta,\eps},\zz)=0$  we have that
$c_1(\calf)=0$ and therefore $\calf$ belongs to
$H^{0,1}(D_{\delta,\eps})/H^1(D_{\delta,\eps},\zz)$, where
$H^1(D_{\delta,\eps},\zz)$ is imbedded to $H^{0,1}(D_{\delta,\eps})$
as it was explained in  Subsection \ref{EXT.chern}. Let $\omega$ be
a $\dbar$-closed $(0,1)$-from  representing $\calf$. Since
$\calf^{\otimes k}$ is extendable onto $K^2_{\delta,\eps}$ we see
that there exists a $\dbar$-closed $(0,1)$-form $\tau $ on
$K^2_{\delta,\eps}$ such that

\begin{equation}
\eqqno(tors-eq1)
\tau|_{D_{\delta,\eps}} = k\cdot \omega + \dbar f + e,
\end{equation}
where $f$ is a smooth function on $D_{\delta,\eps}$ and $e\in
H^1(D_{\delta,\eps},\zz)$. Dividing \eqqref(tors-eq1) by $k$ we
obtain

\begin{equation}
\eqqno(tors-eq2)
1/k\cdot\tau|_{D_{\delta,\eps}} =  \omega + \dbar (1/k\cdot f )+ 1/k\cdot e,
\end{equation}
\ie we see that $\calf\otimes \cale$ extends to $K^2_{\delta,\eps}$,
where $\cale$ is the torsion bundle which corresponds to the
cohomology class $1/k\cdot e\in H^{0,1}(D_{\delta,\eps})$. I.e.,
$\calf\otimes \cale$ is trivial and, therefore, $\calf$ is
isomorphic to the torsion bundle $\cale^{-1}$. But $\cale^{-1}$ can
be only one of $\calf_k(l)$ by construction.

\smallskip\qed

\smallskip To finish the proof of Theorem \ref{tors-ext-thm} all
what is left is to remark that the restriction of $\calf$ to
$D_{\delta , \eps}$ is a torsion bundle, because $\calf^{\otimes k}$
is isomorphic to $\caln$ and the last is trivial on $K^2_{\delta ,
\eps}$. If $\calf|_{D_{\delta , \eps}}$ is trivial itself then
$\calf$ extends to a neighborhood of $\c$. Otherwise it is
isomorphic to one of $\calf_k(l)$. Theorem \ref{tors-ext-thm} is
proved.

\smallskip\qed

\begin{rema} \rm
\label{inv-l}
We proved that a $k$-th root $\calf$ of a globally
defined holomorphic line bundle near a pseudoconcave boundary point
of index $2$ in dimension $2$ possesses a numerical invariant $1\le
l\le k-1$. This $l$ is correctly defined. Indeed if $\calf$ is
defined say on $Y_{\tau}$ for some $\tau >0$ then its extension to
$Y_0$ is determined uniquely, because morphisms of holomorphic
bundles are analytic objects.
\end{rema}

\newsect[sect.CNTR]{Extension of roots across the contractible analytic sets.}

In this section we are going to prove  Theorem \ref{exc-ext-thm}
from the Introduction and derive corollaries from it.

\newprg[prgCNTR.cover]{Extension of analytic coverings}

For the proof of Theorem \ref{exc-ext-thm} we shall need the
removability of codimension two singularities of analytic coverings.
The needed statement is implicitly contained in \cite{NS} and
\cite{De}, but neither of this papers doesn't contain an explicit
formulation. Therefore we shall sketch  the proof of the forthcoming
Theorem \ref{cov-ext-thm}, refereing step by step to the quoted
texts (more to \cite{De}) and keeping notations close to that ones
in \cite{De}. Recall that an analytic  covering is a proper,
surjective holomorphic map $\c:\widetilde{X}\to X$ between normal
complex spaces of equal dimensions.

\begin{thm}
\label{cov-ext-thm} Let $X$ be a normal complex space and $A$ is an
analytic subset of $X$ of codimension at least two. Then every
analytic covering $\widetilde{X}$ over $X\setminus A$ extends to an
analytic covering over the whole of $X$.
\end{thm}

\begin{rema} \rm
\label{non-an-ob} Let us remark that analytic coverings are not
analytic objects. They don't obey neither uniqueness property (P1)
nor the Hartogs type extension property (P2).
\end{rema}
\proof First of all let us make a few reductions. The problem is
local and therefore, since every analytic space can be realized as
an analytic  covering of a domain in a complex linear space, one
reduces the general statement to the case when $X$ is a pseudoconvex
domain in $\cc^n$. After that let us remark that the branch locus
$Y$ of the covering, being a hypersurface, extends through the
codimension two set $A$ by Remmert-Stein theorem. Denote this
extended hypersurface as $\tilde Y$. If $A$ contains points from
$X\setminus \tilde Y$ then, since the complement to $\tilde Y$ in
the smooth locus of $X$ is locally simply connected the
extendability of the covering over such points follows. Therefore we
can suppose that $A\subset \tilde Y$ in the sequel.

\smallskip Let us remark furthermore that all we need is to separate
the sheets of the covering $\widetilde{X}\to X\setminus A$ by
holomorphic functions. I.e., it is sufficient for our goal for every
point $z_0\in X\setminus \tilde Y$ and its preimages $z_1, ..., z_d$
to find a weakly holomorphic function $h\in \calo'(\widetilde{X})$
such that $h(z_1)=1$ and $h(z_2)=...=h(z_d)=0$. Indeed, after that
our problem will be reduced to the extension of an appropriate
symmetric polynomials with holomorphic coefficients from $X\setminus
A$ to $X$, and this is standard. Construction of such $h$  is
exactly what \cite{De} is about. One only needs to remark that the
proof there "doesn't see" a codimension two "hole" $A$ in the branch
locus $\tilde Y$ of the covering.

\smallskip From that place we closely follow \cite{De}. Denote by $\tilde Y'$
the singular locus of $\tilde Y$ and by $\widetilde{X}'$ the part of
the covering situated over $\tilde Y'$. The proof then uses the fact
that $\tilde Y'$ is of codimension at least two in $X$. But then one
can {\slsf include} $A$ to $\tilde Y'$ and follow  the proof of
\cite{De} with these data.

\smallskip Denote by $\widetilde{X}_0$ the unbranched portion of
$\widetilde{X}$, \ie the part of $\widetilde{X}$ situated over
$X\setminus Y$. By $d\nu \deff \c^*(\frac{1}{2})^ndz\wedge d\bar z$
denote the natural volume form on $\widetilde{X}_0$. As the firs
step one finds a holomorphic function $h$ on $\widetilde{X}_0$,
which separates preimages of $z_0$ as above and is square integrable
with respect to $d\nu$ on $\widetilde{X}_0$, see Proposition 1.5 in
\cite{De}. Second, one takes the equation $g$ of he branching locus
of $\widetilde{X}$ and extends the product  $gh$ through
$\c^{-1}\big(\tilde Y\setminus (\tilde Y'\cup A)\big)$ as in Lemma
1.5 from \cite{NS}. And finally one extends the function $h'\deff
gh$ through the points of $\widetilde{X}'\setminus \c^{-1}(A)$ as in
Lemma 1.10 of \cite{De}. Function $h'$ possesses the required
properties.

\smallskip\qed

\begin{rema} \rm
Remark that he branch locus $\tilde Y$ of the extended covering is
the closure of the branched locus of the initial covering, \ie no
new ramifications occur.
\end{rema}

\newprg[CNTR.ext]{Proof of Theorem \ref{exc-ext-thm}.}
As the first step of the proof we state the following:

\begin{lem}
\label{space-ext} Let $X$ be a normal complex space and $A$ a
codimension two analytic subset of $X$. Let a holomorphic line
bundle $\caln$ over $X$ be given. Suppose that $\caln|_{X\setminus
A}$ admits a $k$-th root $\calf$. Then $\calf$ extends as a $k$-th
root of $\caln$ onto the whole of $X$.
\end{lem}
\proof Denote, as usual by $F$ and $N$ the total spaces of $\calf$
and $\caln$ correspondingly. Projections onto the base in both cases
will be denoted as $\pi$. There exists a natural $k$-linear
holomorphic map $p:F\to N|_{X\setminus A}$ defined as
\[
p : (x,\v) \mapsto (x, \v^{\otimes k}).
\]
Mapping $p$ is an analytic covering of order $k$ with branch set
$X\setminus A$. By Theorem \ref{cov-ext-thm} it extends to an
analytic covering $\tilde p:\tilde F \to N$. The branched set of
$\tilde p$ is $X$, projection $\pi : F\to X\setminus A$ obviously
extends to projection $\tilde\pi : \tilde F \to X$ as well as vector
bundle operations. All this follows from normality of the covering
space $\tilde F$ and the fact that $\tilde F \setminus F$ is of
codimension two.

\smallskip The local triviality of $(\tilde F, \tilde\pi , X)$ one
needs to prove over a point $a\in A$ only. Take some $(a, \v)\in
\tilde\pi^{-1}(A)$, $\v\not=0$ and let $(a, \w)\in N$ be its image
under $\tilde p$. Remark that $w\not=0$, otherwise the fiber $\tilde
F_a$ would be contracted by $\tilde p$ to a point and this violates
the properness of $\tilde p$. Since $\tilde p$ is a covering  with
the ramification set $X$ we see that $\tilde p$ is biholomorphic
between some neighborhoods $V\ni (a,\v)$ and $W\ni (a,\w)$. If
$\Sigma \subset W$ is a graph of a local section $\sigma$ of $\caln$
then $\tilde p|_V^{-1}(\Sigma)$ will be the graph of a local section
of $\tilde\calf$.

\smallskip Therefore $\tilde F$ is a total space of a holomorphic
line bundle $\tilde\calf$ extending $\calf$. It is a $k$-th root of
$\caln$ because $\tilde p$ was extended to.

\smallskip\qed

\medskip\noindent{\slsf Proof of Theorem \ref{exc-ext-thm}.}   Let
$\C : X\to Y$ be the contraction. Lemma \ref{space-ext} is now
applicable to the direct images of $\C_*\calf|_{X\setminus E}$ and
$\C_*\caln|_{X\setminus E}$. Both are holomorphic line bundles on
$Y\setminus A$. The second one extends onto $Y$ as a coherent
analytic sheaf by the Theorem about direct image sheaves of Grauert.
Therefore it extends onto $Y$ as a holomorphic line bundle. Denote
this extension (with some abuse of notation) as $\C_*\caln$. By
Lemma \ref{space-ext} $\C_*\calf$ extends as a holomorphic line
bundle onto $Y$ and stays there to be a $k$-th root of $\C_*\caln$.
Denote by $\C_*\calf$ this extension. Now $\C^*\C_*\calf$ will be an
extension of $\calf$ to $X$, which is a $k$-th root of
$\C^*\C_*\caln$, and the last is an extension of $\caln|_{X\setminus
E}$ onto $X$. Theorem is proved.

\smallskip\qed

\newsect[LOG]{Logarithms of line bundles and pseudoconcave domains.}

\newprg[LOG.log]{Extension of logarithms}

Results of this paper in some sense can be reformulated also as
extension of logarithms of holomorphic line bundles. A {\slsf
logarithm} of a holomorphic line bundle $\caln\in H^1(U,\calo^*)$ is
a cohomology class $\calh \in H^1(U,\calo)$ such that $\exp (\calh)
= \caln$. Here $\exp$ stands for the map from the exponential
sequence. From this sequence it is clear that $\caln$ admits a
logarithm on $U$ if and only if $c_1(\caln)=0$, \ie if and only if
$\caln$ is topologically trivial on $U$. This shows that existence
of $\log\caln$ is equivalent to the existence of all roots from
$\caln$.

\smallskip Indeed, consider the exact sequence
\begin{equation}
\eqqno(exact2) 0\to Pic^0(U)\to Pic(U)\buildrel c_1 \over
\longrightarrow NS(U)\to 0.
\end{equation}
The N\'eron-Severi group $NS(U)$ is a discrete subgroup of
$H^2(U,\zz)$, in fact it is equal to $H^{1,1}(U,\zz)$ $\deff
H^2(U,\zz)\cap H^{1,1}(U,\rr)$, see \cite{Che}. Therefore, would
$c_1(\caln)\not=0$ then $\caln$ couldn't admit roots of infinite
number of degrees. Therefore $c_1(\caln)$ should be zero and then
$\caln\in Pic^0(U)$, \ie is topologically trivial. And vice versa,
if $\log\caln$ exists then $\calf_k\deff\exp (1/k\cdot\log\caln)$
will be a $k$-th root of $\caln$.

\begin{nncorol}
\label{log-ext} In the notations of Theorems \ref{root-ext-thm} and
\ref{tors-ext-thm} let $\calh$ be a logarithm of $\caln$ on
$U_{t_0}^+$.

\smallskip\sli If $n\ge 3$, or $n=2$ and $\ind\rho|_{\c}\not= 2$
then $\calh$ extends as a cohomology class to a neighborhood of
$\Sigma_{t_0}$ and stays there to be a logarithm of $\caln$.

\smallskip\slii If $n=\ind\rho|_{\c} = 2$ and $\calh$ doesn't
extend to a neighborhood of $\c$ then there exists a neighborhood
$V$ of $\c$ biholomorphic to $K^2$ such that $\calh|_{K^2\cap
U_{t_0}^+}$ extends to $D$ and as element of $H^1(D,\calo)$ is equal
to $\j (rA)$ for some $r\in \zz$, $r\not=0$.
\end{nncorol}
\proof We use here notations from Introduction. What concerns (\sli
all what is needed to be remarked is that all neighborhoods to which
we extended roots of $\caln$ do not depend on a degree of a root.
This implies the topological triviality of $\caln|_{U_t^+}$ for some
$t<t_0$ close enough to $t_0$. Let $\calh_1$ be some logarithm of
$\caln$ on $\caln|_{U_t^+}$. Then $\calh_1|_{U_{t_0}^+} =\calh +
\j(\z)$ for some $\z\in H^1(U_{t_0}^+,\zz)$. But along the proof we
saw that in this case the natural map $H_1(U_{t_0}^+,\zz) \to
H_1(U_{t}^+,\zz)$ is surjective. Therefore every cohomology class
$\z\in H^1(U_{t_0}^+,\zz)$ is a restriction of some cohomology class
$\tilde\z\in H^1(U_{t}^+,\zz)$. All what is left is to correct
$\calh_1$ by subtracting  from it $\j (\tilde\z)$.

\smallskip Consider the case (\slii . Let $V$ be a neighborhood of
$\c =0$ biholomorphic to $K^2$. $\calh$ extends to a logarithm of
$\caln$ on $U_{t_0}^+\cup D$ by part (\sli . Since $\caln|_D$ is
trivial we see from \eqqref(jakobi2) that $\calh|_D\in
\j(H^1(D,\zz))$, \ie $\calh|_D = \j(rA)$ for some $r\in \zz$.
Suppose $r=0$. Let $\omega $ be a $(0,1)$-form in $U_{t_0}^+\cup D$
representing $\calh$ via Dolbeault isomorphism. Since $\calh|_D=0$
we have that $\omega|_D = \dbar h$ for some function $h$ in $D$.
Take a test function $\phi$ with support in $K^2$ equal to $\adyn$
on $\frac{1}{2}K^2$. Then $\omega_1 \deff \omega - \dbar (\phi h)$
will still represent $\calh$ in $H^{0,1}(U_{t_0}^+\cup D)$. But this
form is smooth on $H^{0,1}(U_{t_0}^+\cup \frac{1}{2}K^2)$. This
means that $\calh$ is extended to a neighborhood of $\c$. Corollary
is proved.

\smallskip\qed

\newprg[LOG.nemir]{Non-extendability of roots from pseudoconcave domains}
In this subsection we shall give the example of Nemirovski. Recall
that a domain $U$ in a complex manifold $X$ is called pseudoconcave
in the sense of Andreotti or, simply {\slsf pseudoconcave} if $U$
admits a proper exhaustion function $\rho :U\to [r_0,+\infty)$ which
is strictly plurisuperharmonic on $U_{r_1}^+\deff \{\rho (x) >r_1\}$
for some $r_1\ge r_0$. A typical example of such $U$ is a complement
to a Stein compact in a compact complex manifold.

\begin{exmp} \rm
Let $S$ be an imbedded surface in $\pp^2$ with a basis of Stein
neighborhoods and such that it is homological to a complex line. Let
$D$ be a Stein neighborhood of $S$.

\begin{rema} \rm
Such $S$ can be obtained by adding to $\pp^1$ three handles and then
perturbing in order to cancel all elliptic points, see \cite{Fo} or
\cite{Ne} for more details. Genus (at least) three comes from the
adjunction inequality:
\[
S^2 + |c_1\cdot S| \le 2g-2,
\]
which is necessary and sufficient condition for an imbedded surface
to have a Stein neighborhood, unless $S$ is a sphere homologous to
zero.
\end{rema}
Then the restriction of $\calo (1)$ to $\pp^2\setminus D$ is
topologically trivial. Indeed, write
\[
H^2_c(D,\zz) \to H^2(\pp^2,\zz) \to H^2(\pp^2\setminus D, \zz).
\]
The first map is a surjection, because the class dual to $S$ goes to
the generator of $H^2(\pp^2,\zz)$. Therefore $H^2(\pp^2,\zz) $ is
mapped to zero by the second arrow. In particular $c_1(\calo
(1)|_{\pp^2\setminus D}) = 0$. Therefore $\calo (1)|_{\pp^2\setminus
\bar D}$ admits all roots and a logarithm. But neither of them can
be extended to $\pp^2$.

\end{exmp}

\newsect[sect.APP]{Appendix. Hartogs domain in a neighborhood of a critical
point.}

We shall explain in this Appendix how to place a Hartogs figure to a
upper level set of a strictly plurisubharmonic function in a
neighborhood of a critical point of index two in such a way that the
associated bidisc is imbedded, contains a neighborhood of this
critical point and its intersection with the upper level set is
connected. After that in the Remark \ref{shadow} we shall explain
once more what happens with "quasi-analytic" objects in this case,
like roots of holomorphic line bundles, which satisfy (P2) but fail
to satisfy (P1). We use the notations of Subsection \ref{EXT.crit}.

\begin{lem}
\label{index-2-lem} $Y_0$ contains an imbedded Hartogs figure
$i(H^2_{\eps})$ such that:

\smallskip\sli the associated bidisc $i(\Delta^2)$ is also imbedded;

\smallskip\slii $i(\Delta^2)\cap Y_0$ is connected;

\smallskip\sliii the union $i(\Delta^2)\cup Y_0$ contains a neighborhood of zero.
\end{lem}
\proof Let $A$ be the constant from \eqqref(y-plus-f), set $t_A =
\frac{A^2-1}{7}$. Without loss of generality we may suppose in what
follows that $A$ is close to $1$, say $1<A<2$. Consider the
following $1$-parameter family of complex annuli
\begin{equation}
\eqqno(ann-t) R_{t,0} \deff \{z\in K^2: (z_1+ t)^2 + z_2^2 =
t_A^2\},
\end{equation}
parameterized by a real parameter $0\le t\le 2t_A$. It is more
convenient to work with the real form of \eqqref(ann-t):
\begin{equation}
\eqqno(ann-re)
\begin{cases}
(x_1 + t)^2 + x_2^2 = \norm{y}^2 + t_A^2,\cr (x_1 + t)y_1 + x_2y_2 =
0.
\end{cases}
\end{equation}

\smallskip Let us list the needed properties of these annuli.

\medskip\noindent{\slsf $1_R$.} {\it For every $t\in [0,2t_A]$ the
curve $R_{t,0}$ intersects the boundary $\d K^2$ of the convex brick
$K^2$ by the side $\d K\times \inter K_A$.} Really, from the first
equation of \eqqref(ann-re) we get that for $\norm{x}=1$ one has
\[
\norm{y}^2 \le 1 + |2x_1t| + t^2 - t_A^2\le 1 + 4t_A + 3t_A^2\le 1 +
7t_A < A^2,
\]
because of the choice of $t_A$. As a result we see that:

\smallskip\noindent{\slsf $2_R$.} {\it Boundaries $\d R_{t,0}$ are
entirely contained in $Y_0$ for all $t\in [0,2t_A]$.} This is
because $\d K\times \inter K_A\subset Y_0$.

\smallskip\noindent{\slsf $3_R$.} {\it $R_{0,0}$ is entirely contained in $Y_0$.}
This is because in that case the same equation gives $\norm{y}^2 =
\norm{x}^2 - t_A^2< \norm{x}^2$ for all $x+iy\in R_{0,0}$.

\smallskip\noindent{\slsf $4_R$.} {\it  Finally $R_{t_A,0}\ni 0$.}

\smallskip From $1_R - 4_R$ we conclude that the holomorphic family of annuli
\begin{equation}
\eqqno(ann-c) R_{t,\tau} \deff K\cap \{(z_1+ t + i\tau )^2 + z_2^2 =
t_A\},
\end{equation}
with $t+i\tau$ in a neighborhood of $U$ of $[0,2t_A]$ in $\cc$,
sweeps up a neighborhood of zero in $\cc^2$ and inherits the same
properties, provided $U$ is thin enough.

\smallskip \smallskip Cut each $R_{t,\tau}$ along $\{x_2=0, x_1\le -t\}$ to
get the disc $D_{t, \tau} = R_{t,\tau}\setminus \{x_2=0, x_1\le
-t\}$ with boundary $\d D_{t, \tau}$. We obtain a holomorphic one
parameter family of discs $D_{t, \tau}$ each catted from the annuli
$R_{t,\tau}$. $D_{t, \tau}$ are mutually disjoint because
$R_{t,\tau}$ are.

\begin{figure}[h]
\label{crit-fig}
\includegraphics[width=2.0in]{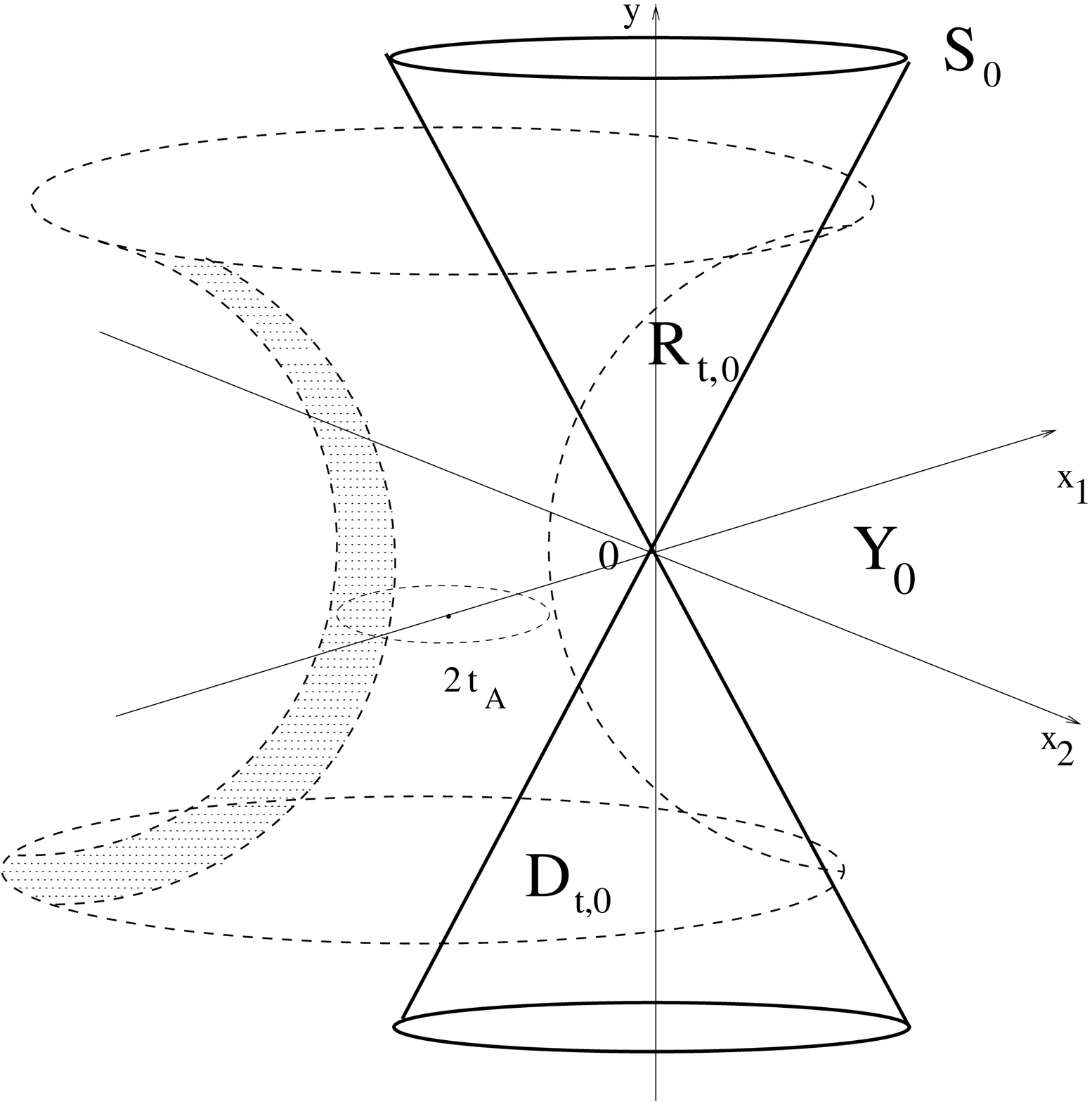}
\qquad\qquad
\includegraphics[width=1.8in]{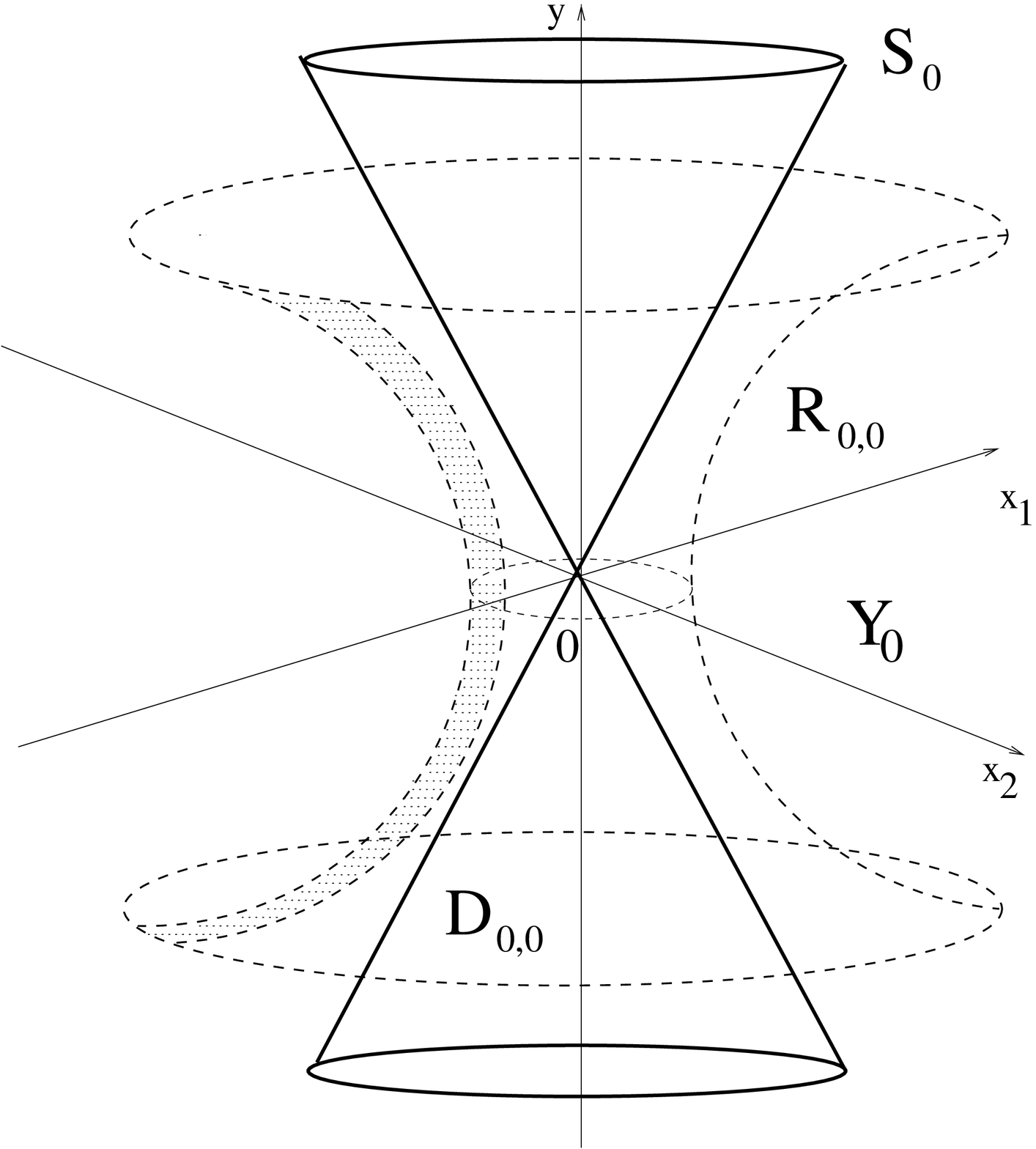}
\caption{On the right the whole annulus $R_{0,0}$ is contained in
$Y_0$ and therefore after cutting $R_{0,0}$, \ie removing of the
dashed zone on the picture, the boundary of the obtained disc
$D_{0,0}$ is still in $Y_0$. When $t\nearrow 2t_A$  the disc
$D_{t,0}$ moves (to the left) sweeping up a neighborhood of zero,
but $\d D_{t,0}$ stays in $Y_0$.}
\end{figure}

We have the following properties of discs $D_{t, \tau}$:

\smallskip\noindent{\slsf $1_D$.} {\it Discs $D_{t, \tau}$  fill in a domain
$0\in D\deff \bigcup_{t+i\tau \in U}D_{t, \tau}$ diffeomorphic to a
bidisc.}

\smallskip\noindent{\slsf $2_D$.} {\it For $t+i\tau $ close to zero
$D_{t, \tau}$ is entirely contained in $Y_0$.}

\smallskip\noindent{\slsf $3_D$.} {\it Boundaries $\d D_{t, \tau}$ are contained in
$Y_0$.}

\smallskip Really, the intersection $R_{t,\tau}\cap \{x_2=0, x_1\le - t\}$ is given by the
equations $(x_1+t)^2 = \norm{y}^2 + t_A^2$ and $(x_1+t)y_1=0$. Since
$x_1=-t$ is impossible we conclude that $y_1=0$ and get
\begin{equation}
\begin{cases}
 y_1 = x_2 =0\cr
(x_1+t)^2 = y_2^2 + t_A^2.
\end{cases}
\end{equation}
Since $x_1$ is negative we get immediately that $y_2^2 < x_1^2$,
\ie that $R_{t,\tau}\cap \{x_2=0, x_1\le t\} \subset Y_0$.

\medskip In what follows we cut $R_{t,\tau}$, from our family, by
$\{|x_2|\le \eps , x_1\le -t\}$ for $\eps >0$ small enough to keep
all properties $1_D-3_D$. The resulting discs we still denote as
$D_{t, \tau}$ and $Y$ stands for $\bigcup_{t+i\tau}D_{t, \tau}$ as
before. Let $V\subset U$ be a sufficiently small disc centered at
zero - such that for every $t+i\tau \in V$ the disc $D_{t, \tau}$ is
entirely contained in $Y_0$. Set

\begin{equation}
\eqqno(h-0) H_0 \deff \bigcup_{t+i\tau \in U}\d D_{t, \tau}\cup
\bigcup_{t+i\tau \in \frac{1}{2}{\bar V}} D_{t, \tau}.
\end{equation}
By $H_{\eps}$ denote the $\eps$-neighborhood of $H_0$ for some small
$\eps >0$.

\smallskip\noindent{\slsf $4_D$.} {\it The intersection $Y\cap Y_0$ is
connected and is contained in the envelope of holomorphy of $H$.}

\smallskip Since the boundaries of $D_{t, \tau}$ stay in $Y_0$ and
$D_{t, \tau}\subset Y_0$ for $t+i\tau \in V$ the only thing to prove
here is the connectivity of $D_{t, \tau}\cap Y_0$ for every
$t+i\tau\in U$. Suppose not. Let $Y_1$ and $Y_2$ be two connected
components of some $D_{t, \tau}\cap Y_0$. Both of them cannot touch
the boundary of this disc, because it is connected ant its
neighborhood is contained in $Y_0$. Therefore, one of them, say
$Y_1$ is relatively compact in $D_{t, \tau}$. The second one could
be supposed to contain the neighborhood of $\d D_{t, \tau}$. Add to
$Y_2$ the catted of band and obtain that $R_{t,\tau}\cap Y_0$ is not
connected. If we set $X_0\deff K^2\setminus Y_0$ then that means
that the compact $R_{t,\tau}\cap X_0$ has bounded components in its
complement in the Riemann surface $R_{t,\tau}$.

\smallskip The domain $Z_0\deff X_0\cap {\it int}(K^2)$ is the domain
of holomorphy and its closure $K^2\setminus Y_0$ is polynomially
convex. Really, one can easily see that curves $\{ R_{t,\tau}:
t>0\}$ pass through any given point in $Y_0$, do not intersect $\bar
X_0$ and leave $K^2$ when $t\nearrow\infty$. Therefore the
polynomial convexity of $K^2\setminus Y_0$ follows from the
Oka-Stolzenberg theorem, see \cite{St}. Intersection of a globally
defined complex curve in $\cc^2$ with the polynomially convex domain
cannot have bounded components in its complement. This is an
immediate consequence of the maximum principle.

\medskip Shrinking discs $D_{t, 0}$ slightly near the boundary (and keeping
$1_D$ - $3_D$) we can suppose that $D_{t, 0}$ is a real analytic
family of closed analytic discs, \ie that there exists a real
analytic map $i:[0,2t_A]\times \bar\Delta \to Y$ such that $i$
holomorphic on second variable and is a biholomorphism of
$\{t\}\times \bar\Delta$ onto $\bar D_{t, 0}$. This statement is
nothing but the Riemann mapping theorem with parameters. The
extension of $i$ to a neighborhood of $[0,2t_A]\times\bar\Delta$
will be the desired map.

\smallskip\qed

\begin{nncorol}
\label{index-2-cor} Let $0$ be a Morse critical point of a strictly
plurisubharmonic function $\rho$ in a neighborhood $U$ of the origin
in $\cc^n$. Set $U^+\deff \{ z\in U: \rho (z) >0\}$. Then $U^+$
contains an imbedded Hartogs figure $i(E^n_{\eps})$ such that:

\smallskip\sli the associated polydisc $i(\Delta^n)$ is also imbedded and
contains the origin;

\slii $i(\Delta^n\setminus E^n_{\eps})\cap U^+$ is connected.
\end{nncorol}

\noindent The proof of this Corollary immediately follows from Lemma
\ref{index-2-lem} after multiplication by sufficiently small
$\Delta^{n-2}_{\eps}$.

\begin{rema} \rm
\label{shadow}
If one deals with an analytic object one can apply condition (P1)
and extend it to the constructed Hartogs figure $i(E^n_{\eps})$.
This will extend it to the neighborhood of the critical point of
$\rho$.

\smallskip Now let us analyze what happens when one tries to extend an
object $\sigma$ in this situation (say a root of a globally existing
bundle) which satisfies (P2) but not (P1). Then
$\sigma|_{i(H^2_{\eps})}$ extends to $\sigma_1$ on $i(\Delta^2)$.
Now we have two objects: $\sigma$ on $Y_0$ and $\sigma_2$ on
$i(\Delta^2)$. They do coincide on $i(H^2_{\eps})$, but why should
they coincide on $i(\Delta^2)\cap Y_0$? And this is required in
order to get an extension of $\sigma$ to a neighborhood of zero.
I.e., the very same problem that was indicated in item (c) of the
Remark \ref{ext-rem} happens.
\end{rema}

\medskip
\begin{rema} \rm
\label{ne-siu}
{\bf (a)} The proof of Lemma \ref{index-2-lem} is
inspired by constructions from \cite{Iv1}, where a qualitative
version of Bochner's tube theorem was proved, \cite{Bo}.

\smallskip\noindent{\bf (b)} Another way to prove this lemma is to
modify the strictly plurisubharmonic function $\rho$ in a
neighborhood of its critical point, as we deed in the proof,
 \ie in such a way that the difference $\bar Y_0\setminus Y_0$ is a totally real disc.
 Then Lemma 2.20 from \cite{Si} will do the job.

\smallskip\noindent{\bf (c)} In both approaches it is crucial to ``catch''
the critical point by {\slsf one} polydisc (with the corresponding
Hartogs figure sitting in $Y_0$). Really, one could try to extend
the bundle directly, as one can extend holomorphic functions in this
situation, see \cite{Iv1} for more details. But then one deals with
the following situation: the ``envelope'' is $\hat H = \Delta\times
A_{1/2,1}$ and the ``Hartogs figure'' is
\begin{equation}
\eqqno(hart-an)
H = \big(\Delta_{\eps}\times A_{1/2,1}\big)\cup
\big(\Delta\times V\big),
\end{equation}
where $V\deff (\d A_{1/2,1})^{\eps}$ is an $\eps$-neighborhood of
the boundary of the annulus $A_{1/2,1}$ in $\cc$, \ie $V$ is a union
of two annuli $V_1, V_2$. And now the proof of the crucial
Hartogs-type Lemma \ref{hart-l} will fail: the domain $U_1=
\Delta_{\eps}\times A_{1/2,1}$ in this case is not simply connected
and one cannot take roots from $F_1$ - we are using the notations of
the proof of that lemma. To my best understanding it is Hartogs
figures of the type \eqqref(hart-an) which appear in Lemma 8.1.1 of
\cite{ST}.
\end{rema}

\ifx\undefined\bysame
\newcommand{\bysame}{\leavevmode\hbox to3em{\hrulefill}\,}
\fi

\def\entry#1#2#3#4\par{\bibitem[#1]{#1}
{\textsc{#2 }}{\sl{#3} }#4\par\vskip2pt}

\end{document}